\documentstyle[12pt]{report} 

\def \dfll {\leaders \hbox to 1em {\hss.\hss}\hfill}
\def\beq {\begin{equation}}
\def\eeq {\end{equation}}
\newtheorem{theorem}{Theorem}
\newtheorem{cor}[theorem]{Corollary}
\newtheorem{rem}[theorem]{Remark}
\newtheorem{defi}{Definition}
\newtheorem{qes}{Question}
\newtheorem{lemma}{Lemma}
\def\blm{\begin{lemma}}
\def\elm{\end{lemma}}
\def\bdf{\begin{defi}}
\def\edf{\end{defi}}
\def\btm{\begin{theorem}}
\def\etm{\end{theorem}}
\def\bpp{\begin{propos}}
\def\epp{\end{propos}}
\def\bQ {\begin{qes}}
\def\eQ {\end{qes}}
\def\btm{\begin{theorem}}
\def\etm{\end{theorem}}
\def\ben{\begin{enumerate}}
\def\een{\end{enumerate}}
\def\bar{\begin{array}}
\def\ear{\end{array}}

\newcommand{\TO}[2]{\stackrel {\mbox{#1}}{\hbox to #2pt{\rightarrowfill}}}

\def\twoheadrightarrow{\to\mkern-20mu\to}
\def\thrafill{$\mathsurround=0pt \mathord- \mkern-6mu 
\cleaders\hbox{$\mkern-2mu
\mathord- \mkern-2mu$}\hfill \mkern-6mu\mathord\twoheadrightarrow$}
\newcommand {\onto} [1]{\hbox to #1pt{\thrafill}}
\newcommand {\ONTO} [2]{\stackrel{\mbox{#1}}{\onto {#2}}}

\newcommand {\prg}[2]{\vspace*{.5cm}\subparagraph{ #1)\ #2\ }\lll{pg-#1}\ \newline}

\def\ub{\underline}
\def \A {\mbox{$\cal A\,$}}
\def \F {I\mkern -6.2mu  F}
\def \V {\mbox{$\cal V\,$}}
\def \C {\mbox{$\cal C\,$}}
\def \Am {\mbox{${\cal A}\!-\!mod\,$}}
\def\id{{1\mkern-5mu {\rm I}}}
\def \li {\mbox{$\cal L\,$}}
\def \lz  {\langle}
\def \rz  {\rangle}

\def \hopf {\mbox{$^0\bigcirc\mkern-13mu\bigcirc^0\,$}}
\def \nhopf {\mbox{$^0\bigcirc\mkern-13mu\bigcirc\,$}}
\def \twhopf {\mbox{$^0\bigcirc\mkern-13mu\bigcirc^1\,$}}

\def\B{{\cal B}}
\def\Bo{{\B}^{\circ}}

\def\rcoa{\leftharpoondown}
\def\lcoa{\rightharpoonup}

\def\Ii {\mbox{\raise .4 ex\hbox{$\int$}$\!\! I$}}

\def \tgl{\mbox{\boldmath${\cal T}\! gl$\unboldmath}}
\def \stA{\mbox{\boldmath$s\!{\cal T}$\unboldmath$({\cal A})$}}

\def \vc {Vect({\bf C})}
\def\lll{\label}

\begin{document} 

\vspace*{1.5cm}

\begin{center}

\section*{Genealogy of Nonperturbative  \\
 Quantum-Invariants of 3-Manifolds: \\ The Surgical Family\footnote
{Expanded version of a talk given at the special session on 3-manifold 
invariants of the conference on Geometry\& Physics in  Aarhus, Denmark, 
 August 1995} }
\bigskip

\medskip

{\large THOMAS KERLER}\\
 \vspace*{1cm}

\bigskip

November, 1995\footnote{Revised Version: December, 1995}

\end{center}
\vspace*{1.7cm}

{\small \noindent{\bf Abstract :} We study the relations between the invariants
$\tau_{RT}$, $\tau_{HKR}$, and $\tau_L$ of Reshetikhin-Turaev, 
Hennings-Kauffman-Radford, and Lyubashenko, respectively. In particular, we 
discuss explicitly how  $\tau_L$ specializes to $\tau_{RT}$ for 
semisimple categories  and  to $\tau_{HKR}$ for Tannakian categories. 
We give arguments for that $\tau_L$ is the most general invariant that 
stems from an extended TQFT.  
We  introduce a canonical, central 
 element, {\sf Q}, for a
quasi-triangular Hopf algebra, $\A$, that allows us to apply the Hennings algorithm
directly, in order to compute 
$\tau_{RT}$, which is  originally obtained
from the semisimple trace-subquotient of $\A-mod$.
Moreover, we generalize Hennings' rules to the context of cobordisms, in 
order to obtain a TQFT for connected surfaces compatible with $\tau_{HKR}\,$.
As an application we show that, for lens spaces and $\A=U_q(sl_2)\,$, the ratio
of $\tau_{HKR}$ and $\tau_{RT}$ is the order of the first homology group.
In the course of this paper we also outline the topology and the algebra that
enter invariance proofs, which contain no reference to 2-handle slides,
but to other moves that are local.
Finally, we give a list of open questions  regarding cellular invariants,
as defined by Turaev-Viro, Kuperberg, and others, their relations among each
other, and their relations to the surgical invariants from above.}
\medskip

\subsection*{Contents}

{\bf 1.) Survey of Surgical Quantum-Invariants}\hfill\pageref{pg-1}

{\em\qquad 1.1) The Reshetikhin-Turaev Invariant }\dfll\pageref{pg-1.1}

{\em\qquad  1.2) The Hennings Invariant }\dfll\pageref{pg-1.2}

{\em\qquad  1.3) The Lyubashenko Invariant }\dfll\pageref{pg-1.3}

{\em\qquad  1.4) Quantum-Invariants from TQFT's }\dfll\pageref{pg-1.4}

{\em\qquad  1.5) $\tau_L$ as a Quantum-Invariant }\dfll\pageref{pg-1.5}

{\em\qquad  1.6) Some Arguments for Exhaustiveness  of $\tau_L$ }
\dfll\pageref{pg-1.6}
\medskip

{\bf 2.) Coends, Universal Liftings, Modularity, Integrals, and All That  }\hfill\pageref{pg-2}

{\em\qquad  2.1) Remarks on Semisimple and Tannakian Categories }
\dfll\pageref{pg-2.1}

{\em\qquad  2.2) General Characterizations of Coends }\dfll\pageref{pg-2.2}

{\em\qquad  2.3) Coends for Special Categories }\dfll\pageref{pg-2.3}

{\em\qquad  2.4) Special Liftings, Pairings, and Modularity }\dfll\pageref{pg-2.4}

{\em\qquad  2.5) Integrals and Cointegrals }\dfll\pageref{pg-2.5}

\medskip

{\bf 3.) 
Lyubashenko's Invariant, and Derivations of $\tau_{RT}$ and $\tau_{HKR}\,$}
\hfill\pageref{pg-3}

{\em\qquad  3.1) Construction of $\tau_L\,$  }
\dfll\pageref{pg-3.1}

{\em\qquad  3.2) Specialization to $\tau_{RT}$ }
\dfll\pageref{pg-3.2}

{\em\qquad  3.3) Categories of Singular Tangles and Hennings' Rules }
\dfll\pageref{pg-3.3}

{\em\qquad  3.4) Computing $\tau_L$ and $\tau_{RT}$ from $\tau_{HKR}$:\newline
\hphantom{Two Fiber F} Two Fiber Functors and a Central Element }
\dfll\pageref{pg-3.4}

{\em\qquad  3.5) Generalization of Hennings' Rules to TQFT's   }
\dfll\pageref{pg-3.5}

\medskip

{\bf 4.)  Open Questions, and Relations with Cellular Invariants}\hfill\pageref{pg-4}

\bigskip

\bigskip

\setcounter{chapter}{1}


\subsection*{1) Survey of  Surgical Quantum-Invariants}\lll{pg-1}

\medskip

At the end of the last decade the studies of quantum-groups,
 subfactors of von Neumann algebras, and the Chern-Simons quantum field theory
brought forward a large and seemingly very powerful new class of so called 
``quantum-invariants'' of closed, compact three-dimensional manifolds. 
In more recent times the focus of research has shifted towards finding 
relations and ``universal'' formulations  for these invariants, which
very well may turn out to be as  instructive and fruitful as 
finding the invariants themselves. A large number of contributions to this 
conference dealt with the connection of the invariant of Reshetikhin-Turaev
with the geometrical data, that we expect from the Chern-Simons functional 
integral, and with other classical invariants by 
using conventional  perturbative methods, or number theoretical expansions.

Yet, there are also still 
plenty  of questions that can be addressed regarding relations among the
invariants that are related to quantum-group constructions. Here, we 
mainly use (non-perturbative) methods from categorical algebra,
and geometric topology. In particular, the considered quantum-groups 
do not have to be deformations of classical Lie algebras, but they do 
have to be finite dimensional. 
\medskip

The main purpose of this contribution is to explain  the relation between 
three ``surgical'' invariants. They are based on surgery presentations of 
three-manifolds, as in [Ki], and use additional algebraic data given 
by, e.g., quantum-groups, or tensor categories, in order to produce a number.
We shall see that the general philosophy of the constructions
is to, first, associate to the complement  $M^{split}\,=\,S^3-U(\li)\,$
of a tubular neighborhood of the surgery link $\li$ a vector in a 
tensor power of a vector space, using only the braiding and rigidity 
information of the algebraic input data. 

The invariant is then 
given by an evaluation against an ``integral'', which corresponds in the surgery operation
to the regluing of the opposite full tori .
The specific assignments are, however, not precisely the  same,
and may lead to different invariants. 
Other types of quantum-invariants, that we shall call here ``cellular
invariants'', start from a presentation of the manifold in terms of 
decompositions into three-dimensional handles. We will come back to
them with a few more details and questions in Chapter 4.
\medskip

In this chapter we summarize the relations between the surgical invariants
$\tau_{RT}$, $\tau_{HKR}$, and $\tau_L$, that will be proven in Chapter 3.
We also give an outline of the general, categorical structure of the axiomatics 
of an extended TQFT, and we will argue that $\tau_L$ is the most general 
invariant that is a specialization of such a structure.  
In the discussion, and also in the rest of the article, 
we shall use the usual notations of tensor categories with  conjugate
objects. The  conventions of the rigidity morphisms are as in [Ke1], which
is the flipped version of [D]. For definitions of the notions of 
abelian categories, tensor categories, etc., see also [McL].
\medskip

\prg{1.1}{The Reshetikhin-Turaev Invariant}
The construction given in [RT]
starts from any  {\em semisimple}, abelian, modular, 
braided tensor category, $\cal C$, with a 
finite set, $\cal J$, of equivalence classes of irreducible objects.
Here, modularity refers to a non-degeneracy of the braiding that we will
define and discuss for general categories in  Paragraph 2.4.

The first step in the construction is to assign to  a generic  planar 
projection of a framed link, $\li\subset S^3\,$, with a generic 
height coordinate, a composite of morphisms, where we associate braid-morphisms
to crossings and rigidity morphisms to maxima and minima.
The objects are specified by an orientation of the link, and a 
coloring, that assigns to the $\nu$-th component an  irreducible object, 
$j_{\nu}\in{\cal J}\,$. For a closed link with $k$ components
we arrive at a morphism $1\to 1$, i.e., a number 
$I_{RT}(\li ;j_1,\ldots,j_k)\,$, which is an invariant of the link for all
colorings. 

A closed compact three-manifold $M$ can  now be presented by a framed link. 
An invariant of $M$ should be, in particular, an invariant of the link. 
As an  Ansatz for $\tau_{RT}\,$ let us therefore use a general linear combination of
the above link invariants
\beq\lll{eq-RT-sum}
\tau_{RT}(M,\C)\;=\;\sum_{j_1,...,j_k\in{\cal J}_0}w(j_1,\ldots,j_k)I_{RT}(\li ;j_1,\ldots,j_k)\qquad,
\eeq
where ${\cal J}_0$ is a system of representatives of $\cal J\,$, and 
we have weights given by a function $w:{{\cal J}_0}^k\to{\bf C}\,$.
It turns out that the additional constraints imposed on the invariant
by Kirby-moves ${\cal O}_1\,$ and ${\cal O}_2\,$ (see [Ki]),
 are solved by a unique  function  
$w$ (for modular \C and up to overall scalings). This will be related  in 
Paragraphs 2.5 and 3.2 to
the existence and uniqueness of integrals of Hopf algebras.
\smallskip

A very common abuse of language in this context is to
associate an invariant, $\tau_{RT}(M,\A)\,$, to a non-semisimple
quantum-group, $\cal A$, as for example
$U_q(\mbox{\tt\bf g})'\,$ with $q$ a root of unity. 

The implied meaning of this is, that we should consider the
category of representations \Am\ and use as input for the 
construction the semisimple trace-quotient 
\beq\lll{eq-mod-quot}
\C_{\A}\,:=\,\overline{(\Am)}^{tr}\;,
\eeq
which has been defined in [Ke1].
\medskip

\prg{1.2} {The Hennings-Invariant}
 Very soon after that M. Hennings 
constructed an invariant, starting from any finite dimensional, modular
 quantum-group,
\A, see [H]. In this approach the representations theory of \A\ is not needed.
Instead a projected link, that describes a manifold, $M$, by surgery,  is mapped 
to an element in a category of singular links, which are decorated with elements
in \A. To this an algorithm is applied, that involves a set of quite intuitive 
combinatorial rules, which had already appeared in [Re]. The elements in \A\ that
are eventually computed are evaluated against a right integral of \A.
 This way we obtain a number, 
$\tau_{HKR}(M,\A)\,$, which turns out to be an 
invariant of $M$. In recent papers, e.g., [KR], L. Kauffman and D. Radford reformulated the 
invariant avoiding the use of orientations and gave a clearer and more
natural proof of invariance under the ${\cal O}_2\,$-move, i.e.,
2-handle-slides. In special examples, with $\A=U_q(s\ell_2)'\,$, they and
T. Ohtsuki [O] established that $\tau_{HKR}\,$ differs from $\tau_{RT}\,$.
Generalization of these results will be
given in Paragraph 3.5 using results from [Ke2]. The following is a
special case of the main result in Theorem~\ref{thm-main-rel} 
stated in  the next paragraph, and proven in Paragraphs 3.2, 3.3, and 3.4.

\blm\lll{lm-RT=H-iffss}\ 

Suppose $\A$ is a finite-dimensional,
modular quantum-group with modular trace-quotient, $\C_{\A}$,
of its representation category. Then the identity 
\beq\lll{eq-RT=H-iffss}
\tau_{RT}(M,\A)\;=\;\tau_{HKR}(M,\A)
\eeq
holds for all closed, (oriented), compact three-manifolds $M$ if 
and only if $\A$ is semisimple.
\elm

If \A\ is not semisimple they differ on most manifolds that are not rational 
homology-spheres, since in this case  $\tau_{HKR}(M,\A)\,=\,0\,$, see [O] 
for $U_q(s\ell_2)'\,$ and [Ke3] for general \A.
A question that makes sense in the light of the known examples is,
whether the invariant are equal for integral homology-spheres,
even if \A\ is not semisimple.

\medskip

A particularly interesting case of a semisimple, modular algebra is the 
Drinfel'd double $D(G)\,$ of a finite group. It has been shown in [AC]
that $\tau_{RT}(M)\,=\,\frac 1 {|G|}\bigl|Hom(\pi_1(M),G)\bigr|\,$, which required
a discussion of the representations theory of $D(G)$. With the above correspondence
a more direct proof is  given in [KR] using the rules for  computing
$\tau_{HKR}\,$. 

In the study of Chern-Simons theory with a finite
gauge group, see [F] and [AC], one is also interested in the computation of
$\tau_{RT}\,$ for the quasi-quantum-group $D(G)^{\alpha}\,$, as defined in
[DPR],  for a non-trivial 3-cocycle of $G$. This corresponds to the 
``quantized'' theory. Moreover, $\alpha$ corresponds to the coupling of the closely 
related quantum field theory from [DW], in a way that we shall explain
a little more in Chapter 4. 

As we shall see in Chapter 3, it is undoubtly possible to extend Hennings'
rules to quasi-Hopf algebras with a non-trivial associator $\phi\in \A^{\otimes 3}\,$. 
Computations in the Hennings setting might make it easier to shed some light on open 
questions regarding the field theory associated to $D(G)^{\alpha}\,$.

\prg{1.3}{The Lyubashenko Invariant}
 In several papers, see [L], V. Lyubashenko
constructs in a more abstract fashion  for any modular, abelian, braided tensor category, 
\C, ``with enough limits'' a set of morphisms, which satisfy the relations of 
the usual generators of the mapping class groups. Moreover, he gives a 
construction of an invariant, $\tau_L(M,\C)\,$, of closed three-manifolds.

 The novelty of this approach  is, that it is  neither assumed that \C\  is semisimple nor that
\C\ is the representation category of any Hopf-algebra. Instead,
{\em categorical} Hopf algebras, obtained from so called {\em coends}, will enter
the construction of the invariant. These notions and related techniques will be 
summarized in Chapter 2. The basic idea of defining an invariant will, however,
be the same. Only now the algebraic input is a  generalization of the data for  
both of the previous two invariants. 

The following correspondence, to be proven in Chapter 3, asserts that also 
the invariant $\tau_L\,$ computed from this 
is a generalization of both $\tau_{RT}\,$ and  $\tau_{HKR}\,$.

\btm\lll{thm-main-rel}\

\ben
\item Suppose \C\ is a semisimple, modular, braided, tensor category, with
$|{\cal J}|\,<\,\infty\,$. 

Then
\beq\lll{eq-rel-L=RT}
\tau_{RT}(M,\C)\,\;=\;\tau_L\bigl(M,\C\bigr)\;.
\eeq
\item Suppose \A\ is a modular quantum-group with $dim(\A)\,<\,\infty\,$. 

Then
\beq\lll{eq-rel-L=HKR}
\tau_{HKR}(M,\A)\;=\;\tau_L(M,\Am)\;.
\eeq
\een
\etm

The missing implication of Lemma~\ref{lm-RT=H-iffss} follows immediately, 
once we
know that $\C_{\A}\,=\,\Am\,$ for semisimple \A. This assertion is not 
entirely obvious, since we also have to make sure that the quantum-dimensions 
of the irreducibles  of \Am\ are not zero. It follows from the fact that
the $K_0$- (or fusion- or Grothendieck- ) ring over ${\bf Z}^+\,$ 
for a rigid, fully reducible  category has no proper ideals. However,  irreducibles
with zero quantum-dimension would form such an ideal.

\prg{1.4}{ Quantum-Invariants from TQFT's}
  There is a  justification 
in calling $\tau_{RT}\,$  a quantum-invariant, since (for certain \C ) it is thought to be the 
same as the invariant defined from the Chern-Simons functional integral, see [Wi]. 
But for $\tau_L\,$ and $\tau_{HKR}$, as well as $\tau_{RT}$ obtained from other categories,
similar action-functionals  are not known and, in fact, are not likely to exist.
We recall, however, that the characteristics of what we consider a quantum field theory in
2+1 space-time dimensions may also  be given in the canonical or algebraic formalism.
 
The basics include, that we associate to a  surface, $\Sigma_t$, given by space points at a given 
time, a Hilbert space of states or, with more general axioms, a vector space, $\V(\Sigma_t)\,$.
In the Hamiltonian picture we also have a propagator, i.e., a unitary map,  between the vector spaces
at different times.

Atiyah [A] generalized this axiom for {\em topological} quantum field theories (TQFT's) by
associating such a map, $\V(M)\,:\,\V(\Sigma_1)\to\V(\Sigma_2)\,$ to any three-fold that
cobords the parametrized surfaces $\Sigma_i$ to each other. 
The time-coordinate is now a Morse function.
Dropping the assumptions,
regarding hermitian structures, a TQFT is thus nothing else but a fiber-functor on the 
category $Cob_3\,$ of 2+1-dim cobordisms. If we assume $\V(\emptyset)={\bf C}\,$ then
this assigns in fact a number, $\V(M)\,=\tau(M)1\,$, to each closed manifold, $M\,$. 

The precise definition also involves a central extension of the cobordism category 
(or, equivalently, the use of projective functors),  which we shall suppress in the
following discussion. 
In this sense we shall give the following, provisional definition:

\bdf\lll{defi-Q-inv}
We say that $\tau:M\mapsto\tau(M)\in{\bf C}\,$ is a \ub{quantum-invariant} if it  
specializes from an extended TQFT, $\,\V\,:\,\, Cob_3(*)\to {\bf AbCat}\,$.
\edf

The notion of an extended TQFT used here is motivated by the situation of Chern-Simons
with field insertions, or, in the general language of algebraic field theory, 
the inclusion of other sectors, that may not have vacua. The mathematical axioms
require that we associate to every compact, one-fold, $S$,
 an abelian  category, $\C(S)\,$. This assignment
shall be compatible with the tensor product structures given 
for abelian categories by Deligne's product $\odot$ (see [D]), and for the topological categories by  disjoint unions.
To a surface, $\Sigma\,$, with $\partial \Sigma\cong S\,$, we associate an object, $\V(\Sigma)\,$,
 in $\C(S)\,$, and to cobordisms between these we associated morphisms in $\C(S)\,$. 

Suppose that $Cob_3(N)\,$ is the category of cobordisms between compact surfaces with 
$N$ holes. Then the axioms  require a functor $\V_N\,:\,Cob_3(N)\,\to\,\C(\amalg^NS^1)\,=\,
\C(S^1)\odot\ldots\odot\C(S^1)\,$. The ideas for  formulating an extended TQFT in this way go back 
to D. Kazhdan and N. Reshetikhin, [Ka], and have been put in the  concrete setting of Chern-Simons
theory with finite gauge group (as well as into  writing)  by D. Freed [F]. 

They imply Atiyah's definition if we set $\C(\emptyset)\,=\,Vect({\bf C})\,$.
There is indeed little choice in this  assignment, since the category of 
vector spaces is the only canonical one that acts as a unit with respect to the
$\odot$-multiplication. 
\medskip

If we distinguish on  a surface between start and target holes, there is a canonical way
 to reinterpret the object
$\V(\Sigma)\in \C(S_s\amalg S_t)\,$ as a functor, $\V(\Sigma)\,:\,\C(-S_s)\to \C(S_t)\,$;
see the introduction of  [Ke2] or also [KL].
Similarly, the morphisms may be rewritten as  natural transformation. This explains the
notation in Definition \ref{defi-Q-inv}, where we think of $\V$ as a 2-functor from
the 2-category of 1+1+1-dimensional cobordisms to the 2-category of abelian categories.

Although this language might be somewhat  unfamiliar to the non-categorist, it is in many
ways a more natural point of view, and we obtain a few more, very useful constraints on $\V\,$.
Aside from a few arguments that we will need in the next paragraph, extended TQFT's will
not be used in the remainder of this paper. We will, however, continue to discuss ordinary
TQFT's, e.g., in Paragraph 3.5.

\prg{1.5}{ $\tau_L$ as a Quantum-Invariant}
An invariant that is given by the  (renormalized) partition function of a functional 
integral, such as Chern-Simons, is quite obviously expected to be a quantum-invariant.
Recall also, that if $\tau(-M)\,=\,\overline{\tau(M)}\,$, we can often (re-)construct
a TQFT from the invariant alone by an obvious generalization of
the GNS-procedure. (The signed, inner product spaces, that are first obtained, can
be quite huge and need not be self-dual so that we usually have to impose further
regularity conditions on $\tau$). Moreover, constructions have been proposed by
Turaev in [T], in which a TQFT is associated to $\tau_{RT}$, starting again from a 
semisimple category.

Despite the similarities discussed in the previous sections TQFT's do not exist for $\tau_L\,$
in this na\"\i ve way. In particular we will see at the end of Paragraph 3.3 that
\beq\lll{eq-S1S2=0}
\tau_{HKR}\bigl(S^1\times S^2,\A)\;=\;0\qquad {\rm iff\ \A\ not\ semisimple.}
\eeq
For a TQFT in the usual framework this also implies that $\V(S^2)=0\,$, which
is possible only if $\V\equiv 0\,$. 

Still, it can be shown that the consistency problem arising here is solely one of treating 
the connectivity of cobordisms and surfaces in the right way. The easiest way to
circumvent it is  to consider the categories $Cob_3^{conn}(N)$, which consist 
only of connected cobordisms between connected surfaces. We may of course also consider 
disjoint unions of such cobordisms, but the important features of a tensor category, 
such as rigidity, leading  to the  contradiction in (\ref{eq-S1S2=0}),  are no longer present.
In this restricted setting the following result permits us to view $\tau_L$ still 
as a quantum-invariant:

\btm [{[KL]}]\lll{thm-conn-TQFT} 
For any abelian, rigid, modular, balanced, braided tensor category, \C, with
certain limits, we have a series of functors
\beq\lll{eq-conn-TQFT}
\V_N\;:\;Cob_3(N)^{conn}\,\;\longrightarrow\,\C^{\odot N}\;,
\eeq
which respect the tensor products and 2-categorical compositions.
\etm

In this construction the choice of the circle category $\C\,=\,\C(S^1)\,$ determines
the construction completely. Also, \V\ specializes to the  $\tau_L$, and we reconstruct
the same representations of the mapping class groups as in [L], only now directly 
derived from  tangle-presentations of  the cobordism categories as in [Ke5].

In Paragraph 3.5 we shall give a construction of the functor $\V_0\,$ 
in the case $\C\,=\,\Am\,$. Instead of the categorical picture in [KL],
we  use here a translation into the combinatorial Hennings picture,
extended to manifolds with boundary.
\medskip

Let us conclude this paragraph with a few remarks on the construction of so called
{\em half-projective} TQFT's
in the disconnected, non-semisimple case, given a TQFT for connected surfaces as in [KL]:

As in the semisimple case, we first write a connected cobordism, $M:\Sigma_s\to\Sigma_t\,$,
in the form $M\,=\,\Pi_{\Sigma_t}^{\dagger}\circ\widetilde M\circ \Pi_{\Sigma_s}\,$, where
$\widetilde M$ is a cobordism between connected surfaces. Here we used a choice of 
cobordisms, 
$\Pi_{\Sigma}:\Sigma_1\amalg\ldots\amalg\Sigma_K\,\to\,\Sigma_1\#\ldots\#\Sigma_K\,$,
for every surface, $\Sigma$, with $K$ connected components $\Sigma_j\,$. The cobordism
$\Pi_{\Sigma}^{\dagger}\,$ is an arrow in the opposite direction with the property that 
$\Pi_{\Sigma}^{\dagger}\circ\Pi_{\Sigma}\,=\,\id_{\Sigma_1}\#\ldots\#\id_{\Sigma_K}\,$.
>From this we easily see that $\Lambda_{\Sigma}\,:=\,\Pi_{\Sigma}\circ\Pi_{\Sigma}^{\dagger}\,$
obeys
\beq\lll{eq-Lambda-2}
\Lambda_{\Sigma}\circ\Lambda_{\Sigma}\,=\,
\Lambda_{\Sigma}\#\underbrace{(S^1\times S^2)\#\ldots\#(S^1\times S^2)}_{K-1\,{\rm\ times}}
\qquad\;{\rm and}\;\qquad 
\Lambda_{\Sigma}\in Cob_3^{conn}\;.
\eeq
In the semisimple case $\tau(S^1\times S^2)$ can be normalized to 1 so that a TQFT, $\cal V$,
associates to
$\Lambda_{\Sigma}$ a projector. The map ${\cal V}(M)\,$ is the reduction of 
${\cal V}(\widetilde M)\,$ to the respective subspaces 
$V_{\Sigma_{s/t}}:=im\bigl({\cal V}(\Lambda_{\Sigma{s/t}})\bigr)\,$.

In the non-semisimple case it follows for $K\geq 2$ that ${\cal V}(\Lambda_{\Sigma})\,$
is merely nilpotent of order two, and in fact will be zero for a slight variation of the
connected TQFT from [KL]. In [Ke3] we use this observation to construct a map 
${\cal V}:Cob_3\to \vc\,$, which extends the functor on $Cob_3^{conn}$ and has all the 
properties of a functor, except that it is has half-projective compositions. This means
\beq\lll{eq-half-ext}
\V\bigl(M_2\circ M_1\bigr)\,\;=\;\,{\tt x}^{\mu(M_2,M_1)}\,\V(M_2)\V(M_1)\qquad,
\eeq
where $\mu$ is a cocycle on $Cob_3$ with values in ${\bf Z}^{0,+}\,$, and {\tt x} is
not necessarily invertible. In our case we have ${\tt x}=0$ (recall that $0^0=1$), 
and $\mu$ is $K-1$ if two connected cobordisms are glued  over $K$ boundary components. 

These properties also imply that a cobordisms, whose ``interior'' homology or ``interior''
fundamental group (see [Ke3]) is non-trivial, is mapped to zero. This generalizes a vanishing 
result for closed manifolds and $U_q(s\ell_2)\,$, that was obtained in [O] by 
a more or less direct computation.

In the following, if we speak of extended TQFT's or quantum-invariants, we shall assume
that we have taken care of connectivity questions in some way. In particular, we count $\tau_L\,$
as a quantum-invariant.

\prg{1.6}{ Some Arguments for Exhaustiveness of $\tau_L\,$}
Theorem~\ref{thm-main-rel}  does in fact suggest that $\tau_L\,$ is a most general 
invariant in some class of invariants. Looking at Definition~\ref{defi-Q-inv} it is
a-priori not clear why this should be the class of quantum-invariants. 
For example the definition does not refer to any braided tensor structure or 
quantum-groups. One might thus hope to find TQFT's coming from more general and
more exotic categories.

This is, however, not possible.  The mentioned  quantum-algebraic structures are already  
contained in the  topology of  three-dimensional manifolds with corners.
In particular, the generating category $\C=\C(S^1)\,$ must reflect the
properties of $Cob_3(1)\,$. In [Ke4] we will give an abstract algebraic
characterization of this topological category as follows:

\btm\lll{thm-alg=cob} Suppose $\cal G$ is \ub{the} free strict, balanced, braided 
tensor category, freely generated by a self-dual, braided 
Hopf algebra-object, $\ub F$,
with integrals and a non-degenerate Hopf pairing, $\ub F\otimes \ub F\to 1\,$. 

Then there is a surjective functor:
$$
{\sf T}\;:\;{\cal G}\,\ONTO{}{20}\,Cob_3(1)^{conn}\;.
$$
(Here the topological category shall have only one object in each 
isomorphism class.)
\etm

Since elementary Hopf algebra relations are quite naturally identified 
by {\sf T} with elementary Cerf theoretical relations, we are led to
 conjecture here that {\sf T} is in fact an {\em isomorphism} of categories.
\medskip

Hopf algebra structures - in a topological sense -  
for the punctured torus were  independently also formulated by Crane and 
Yetter in [CY]. The generators described there, however, are only those
appearing in the axioms for Hopf algebras. Consequently, the corresponding 
cobordisms are  only those that are embeddable into ${\bf R}^3\,$.
Yet, in order to guarantee the asserted surjectivity, we also need to 
include  generators,
which account for  surgery  on cobordisms. Those are precisely the integrals,
entering the definition of $\cal G\,$. Moreover, the Hopf pairing is used
to describe a Dehn twists between neighboring one-handles. 
It is closely related
to the modularity requirement, to be discussed in Section 2.4. 
For these reasons the proof of Theorem~\ref{thm-alg=cob} cannot be done 
using simply diagrams on a surface, as in [CY], but higher dimensional 
presentations,  as in [Ke5], have to be employed. 
\medskip

It is relatively easy to see that gluing surfaces to a three holed sphere,
$\Sigma_{0,3}$, 
allows us to define a tensor product in $Cob_3(1)^{conn}$.
There are also
easy choices of cobordisms that make it into a strict and braided category,
see, e.g.,  [KL].

The generating object is  obviously the punctured torus,
${\sf T}(\ub F)\,=\,\Sigma_{1,1}\,$. Since the classical notion of a 
Hopf algebra makes sense  only in  a symmetric category, we have to 
use the natural, braided modification introduced by Majid [M1]. 
The generating morphisms in $\cal G$, like (co-)products and (co-)integrals, 
can be naturally represented by elementary cobordisms, i.e., ones that
correspond to simple handle-attachments of dimensions three or four. 
\medskip

Now, the composite $\V_1\circ{\sf T}\,:\,{\cal G}\,\to\,\C\,$ shows that
$\C\,$ must have all the properties required in Theorem~\ref{thm-conn-TQFT}.
For example we have a tensor structure given by
\beq\lll{eq-def-tensor}
\otimes\,:=\,\V(\Sigma_{0,3})\;:\C\odot\C\,\longrightarrow\,\C\qquad.
\eeq

The only thing  that still needs explanation, is the  relation between
 the existence of a Hopf algebra object and 
the existence of certain limits.  Specifically, the limit we need is going
to be the following {\em coend}:
\beq\lll{eq-coend-not}
F\;:=\;\int X^{\vee}\otimes X\quad\in\,\C\qquad.
\eeq
 
We shall give a meaning of to this formula in Chapter 2.  Let us 
outline here, how this formula follows from  the axioms of an extended TQFT
and basic topology. We shall freely use
the coend notation; the details of its definition are not relevant
to the general line of arguments presented next:
\medskip

The disc, $D^2$, seen as a 1+1-cobordism from
$S^1$ to $\emptyset$ induces by the axioms a basic fiber-functor
${\cal I}\,:=\,\V_1(\Sigma_{0,1})\,:\,\C\,\to\,Vect({\bf C})\,$.
The choice that appears to be the most natural one coming from
the physical examples  is the invariance-functor:
\beq\lll{eq-I=Inv}
{\cal I}=Inv\,:\; X\,\mapsto\,Hom_{\cal C}(1,X)
\eeq

Although there are some constraints on $\cal I\,$, this choice is by 
no means mathematically stringent. A slightly different fiber-functor is 
for example 
chosen in [Ke3]. Still, let us assume it here as an additional axiom of 
extended TQFT's.

Now, the two-holed sphere can be thought of as a cobordism 
$\Sigma_{0,2}\,:\,S^1\amalg S^1\,\to\,\emptyset\,$, which in turn is 
the composite of the pair of pants $\Sigma_{0,3}$ and a disc. Using the
representations of the latter in (\ref{eq-def-tensor}) and (\ref{eq-I=Inv}).
we find that $\Sigma_{0,2}$ is represented by the following functor:
\beq\lll{eq-Omega}
\Omega\,:\;\C\odot\C\,\longrightarrow\,Vect({\bf C})\,\,:\;X\odot Y\,\mapsto\,
Inv(X\otimes Y)\;.
\eeq

But  $\Sigma_{0,2}$ has an obvious, two-sided
 inverse given by the same cylinder 
in opposite direction,  $-\Sigma_{0,2}\,:\,\emptyset\,\to\,S^1\amalg S^1\,$.
The existence of an extended TQFT thus implies that $\Omega$ has a two-sided
inverse, $\Omega^{-1}\,:\,Vect({\bf C})\,\to\,\C\odot\C\,$. A functor of 
this form is up to isomorphisms given by an object $\F\,\in\,\C\odot\C\,$
(e.g., the image of $\bf C$). It turns out that $\F$ inverts $\Omega$ iff
it is the coend:
\beq\lll{eq-FF}
\F\;=\;\int X^{\vee}\odot X
\eeq

This does not always  exist (like $\Omega$ may not be invertible 
in some categories). In order to say it is an inverse we actually also have to specify 
an isomorphism of the composite to the identity, which gives rise to the transformations 
that are also part of the definition of a coend.
The precise meaning of the additional 
requirement about ``enough limits'' is thus that $\C$ contains the element
$\F\,$. 
We shall see in the examples of Paragraph 2.3 that it has to be understood as a finiteness
condition for $\C\,$.
\medskip

In light of Theorem \ref{thm-alg=cob} we should be really interested in the 
functor associated to $\Sigma_{1,1}\,:\, \emptyset\to S^1\,$, which is given
by an object $F$  in  $\C\,$. We notice that $\Sigma_{1,1}\,$ is the composite
of $-\Sigma_{0,2}\,$ with $\Sigma_{0,3}\,$ over $S^1\amalg S^1$. This results 
in the following identity, from which we also obtain formula 
(\ref{eq-coend-not}):
\beq\lll{eq-F=tFF}
F\;=\;\otimes(\F)\qquad.
\eeq

At this point one should compare the statements of 
Theorem~\ref{thm-conn-TQFT} and Theorem~\ref{thm-alg=cob}, 
keeping in mind that, so far, we have made no assumptions on the structure of
$F$ in  $\C\,$. The observant reader will notice that
this must imply a non-trivial theorem on the structure of coends in abelian,
braided tensor categories in general. The corresponding, more precise 
statement is given in the next theorem, and enters crucially in the construction
of $\tau_L\,$:

\btm[{[L]}]\lll{thm-coend=alg}
Suppose that $\C$ is an abelian, rigid, balanced, modular,
braided tensor category,
for which the coend $F\,$, as in (\ref{eq-coend-not}), exists.

Then $F$ has the structure of a braided Hopf-algebra in $\C\,$.
It has   unique, categorical \linebreak (co-~)integrals, and its non-degenerate 
 pairing is a pairing of Hopf-algebras. In short, there is a functor 
$$
{\sf R}\,:\;{\cal G}\,\to\C\qquad,
$$
which maps $\ub F$ to $F\,$. 
\etm

It is intriguing to see that, although we were 
starting from quite different and very general assumptions, we
can extract, both in the topological as well as in the abelian setting,
distinguished objects, $\Sigma_{1,1}$ and $F$, respectively, which
automatically have the properties of $\ub F$. This may be seen as
the abstract reason, why it is possible to construct the functor from 
Theorem~\ref{thm-conn-TQFT}.
\medskip

It is not hard to see from the constructions, that the 
 functors {\sf T} and {\sf R} are also unique up to obvious isomorphisms
and a few choices of orientations. It thus follows that the extended
TQFT  $\V\,$ is basically unique, too. 

The line of arguments given in this paragraph leaves thus no alternatives 
to the construction of quantum-invariant as in 
Definition~\ref{defi-Q-inv}. It should now be a matter of filling
in the details and making the proper formalization that would 
complete them into a proof of the uniqueness of~$\tau_L\,$.
\bigskip

\setcounter{chapter}{2}

\subsection*{2) Coends, Universal Liftings, Modularity, Integrals, 
and All That}\lll{pg-2}
\medskip

In this chapter we will discuss  some of the properties 
of a coend as in (\ref{eq-coend-not}), especially the consequences of its
universality, which are crucial in the construction of $\tau_L\,$.
We will compute $F$ explicitly for semisimple categories and for 
representations categories. Moreover, we shall determine the 
integrals and cointegrals of $F$ in these cases.
Let us start with some general, preparatory remarks on the 
the two types of categories in question:

\prg{2.1}{ Remarks on Semisimple and Tannakian  Categories\ }
An important fact about a  rigid, balanced, braided tensor category, \C, is that
it always admits a canonical system of traces,
\beq\lll{eq-tr}
tr_X\;:\;End_{\cal C}(X)\,\longrightarrow\,{\bf C}\;,
\eeq

which is generally cyclic and  respects the $\otimes$-product. 
Important ingredients in the construction are the usual rigidity morphisms,
$ev_X\,:\,X^{\vee}\otimes X\to 1\,$ and $coev_X\,:\,1\to X\otimes X^{\vee}\,$,
but also a  canonical set of {\em flipped} rigidities,
\beq\lll{eq-flip-rig}
\widetilde{ev}_X\,:\,X\otimes X^{\vee}\,\to\, 1\qquad{\rm and}\qquad\; 
\widetilde{coev}_X\,:\,1\,\to\, X^{\vee}\otimes X\quad,
\eeq
that are constructed from the ordinary ones using the balancing and braiding 
in \C,  see, e.g.,  [Ke1].
The trace allows us to define the so called quantum-dimensions
\beq\lll{eq-q-dim}
d(X)\;=\;tr_X(\id_X)\qquad\;{\rm for\ all\ objects}\;X\quad,
\eeq
with the properties $d(X\oplus Y)=d(X)+d(Y)$ and $d(X\otimes Y)=d(X)\cdot d(Y)\,$.
Also, we obtain a symmetric pairing
\beq\lll{eq-hom-pair}
Hom_{\cal C}(Y,X)\otimes_{\bf C}Hom_{\cal C}(X,Y)\,
\stackrel {\circ}{\longrightarrow}\,End_{\cal C}(X)\,\stackrel {tr_X}
{\longrightarrow}\,{\bf C}\qquad.
\eeq

As in [Ke1] we say that a category is {\em semisimple} iff this pairing is 
non-degenerate for all $X$ and $Y$. It follows that a morphism, $I:X\to Y\,$,
is completely determined by the matrix elements $gIf\in{\bf C}\,$, where  
$f\in Hom_{\cal C}(j,X)\,$, $g\in Hom_{\cal C}(X,j)\,$, and $j$ is irreducible.
More precisely,  we have the following isomorphisms.
\beq\lll{eq-hom-iso}
Hom_{\cal C}(X,Y)\;\cong\;\bigoplus_{j\in {\cal J}_0}Hom_{\bf C}\bigl(Hom_{\cal C}(j,X),
Hom_{\cal C}(j,Y)\bigr)\qquad.
\eeq
Here ${\cal J}_0\,$ is a set of representatives of the set of 
classes of irreducibles $\cal J$ as in Paragraph 1.1. 
\medskip

Another class of categories are the so called Tannakian ones, i.e., those which admit 
 an exact
$\otimes$-functor into the category of vector spaces:
\beq\lll{eq-tann-fct}
{\sf V}\,:\;{\cal C}\,\longrightarrow\,\vc\quad.
\eeq 
The cardinal example for such a category is given by the representation category
\Am\ of  a Hopf algebra. {\sf V} is then simply the  ``forgetful functor''.

An obvious, necessary condition for a category, $\C$,  to be Tannakian, is that we have
dimensions, $X\mapsto d_X$, which assign to each object a value in ${\bf Z}^{+,0}\,$,
and which respect the sums and products of objects in the way the canonical dimensions from above
do. It follows from Perron-Frobenius theory that for categories, with finitely many isomorphism
classes of irreducibles, there exists exactly one such dimension with values in ${\bf R}^{+,0}\,$.
As a matter of fact, most of the categories $\C_{\cal A}\,$, with 
$\A=U_q({\tt g})\,$, admit a dimension that is  positive but not integral, and hence 
these $\C_{\cal A}\,$ can be 
excluded from the list of Tannakian categories. 

In the case of a symmetric category, \C, and if the {\em canonical} dimensions take values in 
${\bf Z}^{+,0}\,$, it is a result of [D] that \C\ is Tannakian.

A very important question is, whether a Tannakian category is always the of the form \Am\ for 
some Hopf-algebra \A. Given the functor {\sf V}, a natural candidate for the algebra is 
$\A=End_{Cat}({\sf V})\,$. Under a few technical assumptions (which, e.g.,  allow us 
identify \Am\ with $\A^*\!-\!comod$) reconstruction results are obtained - in several stages of
generality - in [Ta]. In this sense we shall often 
use the notions of a representation category and
a Tannakian category synonymously.
\medskip

It is quite useful to try to {\em represent} fiber-functors as in (\ref{eq-tann-fct}).
This means we want  to find an object $Q$ in \C\ itself such that $\cal F$ can
also be expressed by the functor $Hom_{\cal C}(Q,-)\,$, which associates to each 
$X$ the vector space $Hom_{\cal C}(Q,X)\,$.

If we have already made the identification $\C=\Am\,$, with $dim(\A)<\infty\,$,
it is easy to see that $Q$ is given, as a vector space, by \A\ itself, and as module 
it is equipped with the left regular action $a.q=aq\,$. 

The dual $Q^{\vee}\,$ is thus, as a vector space, given by the space of functions 
${\cal A}^*\,$. The action is given by $a.\rho\,=\,\rho\rcoa S(a)\,$ for all
$a\in\A$ and $\rho\in{\cal A}^*\,$. Here, $S$ is the antipode, and $\rcoa$ is 
defined by $(\rho\rcoa b)(y)\,=\,\rho(by)\,$.
\medskip

It is not true that $Q$ is an algebra in \C\ since the multiplication does not 
interwine the action of \A. However, since the comultiplication intertwines,
 $Q^{\vee}\,$ is an algebra, which is commutative in symmetric categories. 
This property is crucial in the  construction of a fiber functor as in [D].

Another algebra, in fact even a Hopf algebra, in \C\ was proposed in 
Theorem~\ref{thm-coend=alg}. We shall see that $F$ is also given by ${\cal A}^*\,$,
but carries the coadjoint action instead.

\prg{2.2}{ General Characterizations of Coends\ }
Quite generally, coends are associated to a pair of categories, 
\C\ and $\cal B$, and a  functor
\beq\lll{eq-S-func}
{\sf S}\,:\;\C^{opp}\times \C\,\longrightarrow\,{\cal B}\;,
\eeq
and are denoted by $\int{\sf S}(X,X)\,$. The aim of this paragraph is to make the
elegant but also rather abstract definition in, e.g., [McL] more concrete in the example,
where ${\cal B}=\C\,$ and ${\sf S}(Y,X)=Y^{\vee}\otimes X\,$, and thus make sense  
of the notation in (\ref{eq-coend-not}).
\medskip

We start by introducing the category $dn({\sf S})\,$ of ``dinatural transformations
of {\sf S} to a constant''. Specifically, this means that the objects of
$dn({\sf S})\,$ are given by pairs $\lz Z,\xi\rz\,$. The first argument $Z$ is an
object in $\C\,$, and the second is a map that assigns to every object $X$  in 
$\C\,$ a morphism $\xi_X\,:\,X^{\vee}\otimes X={\sf S}(X,X)\,\to\,Z\,$. They
shall be such that for any pair of objects, $X$ and $Y$, and any morphism 
$f:X\to Y$ the following diagram commutes:
\beq\lll{eq-dinat-diag}
\bar{ccc}
Y^{\vee}\otimes X\,&\stackrel{\id_Y\otimes f}{\hbox to 45pt{\rightarrowfill}}&\;
Y^{\vee}\otimes Y\\
^{f^t\otimes \id_X}\Biggl\downarrow\;&&\qquad\Biggl\downarrow \ ^{\xi_Y}\\
X^{\vee}\otimes X\,&\stackrel{\xi_X}{\hbox to 45pt{\rightarrowfill}}&\;Z\\
\ear
\eeq

A morphism, $h\,:\,\lz Z,\xi\rz\,\to\,\lz Z',\xi'\rz\,$, in $dn({\sf S})$ is first of 
all given by a
morphism $h\in Hom_{\cal C}\bigl(Z,Z'\bigr)\,$. It shall also map the second arguments into 
each other, i.e., we have the following, additional  requirement: 
\beq\lll{eq-dn-mor}
h\circ\xi_X\;=\;\xi_{X'}\qquad\;{\rm for\ all\ objects\ }\quad X\qquad.
\eeq

It is often instructive to consider a slightly different but equivalent definition of
$dn({\sf S})$. The objects here are again pairs $\lz Z,\delta\rz\,$, but now the second
is a  natural transformation, $\delta\,:\, id \,\dot{\to}\, -\otimes Z\,$, i.e., a
family of morphisms $\delta_X\,:\,X\,\to\,X\otimes Z\,$ with 
$\delta_Y\,f\,=\,f\otimes\id_Z\,\delta_X\,$ for $f:X\to Y\,$. Because of the intertwining 
relation we may think of such a transformation as a {\em coaction} of  $Z$ on the objects of
$\C$, if it were also compatible with a coproduct of $Z$.
\medskip

As we already pointed out, a very  important property of a coend for the 
construction 
of $\tau_L\,$ is the {\em universality property}. More precisely, we say that $\lz F,i\rz\,$
is a coend iff it is an {\em initial} object 
 in $dn({\sf S})\,$.  This means that 
 for any other object $\lz Z,\xi\rz\,$ there exists {\em exactly one} morphism
$r^Z\,:\,\lz F,i\rz\,\to\,\lz Z,\xi\rz\,$.

In other words, whenever we have a system of morphisms that fulfill the dinaturality
condition from (\ref{eq-dinat-diag}) there exists a map $r^Z\,:\,F\to Z\,$, such that 
each of them factors through $r^Z\,$ as below.
\beq\lll{eq-xi-fac}
\xi_X\,:\,X^{\vee}\otimes X\,\stackrel{i_X}{\longrightarrow}\,F\,
\stackrel{r^Z}{\longrightarrow}\,Z\qquad.
\eeq

It follows immediately from general nonsense that if $F$ exists it is also unique
up to isomorphisms.
\medskip

The transformations $\,i_X$ allow us define for every object, $X\,$,
 in a balanced,  braided category a canonical 
invariance, $tr^q_X\in Hom_{\cal C}(1,F)\,$,  which is given by
\beq\lll{eq-q-tr}
tr^q_X\,:\;1\,\TO{$\widetilde{coev}_X$}{50}\,X^{\vee}\otimes X\,\TO{$i_X$}{30}\,F\quad.
\eeq
It has properties similar to those of the trace in (\ref{eq-tr}), e.g., 
we have $tr^q_{X\oplus Y}\,=\,tr^q_X+tr^q_Y\,$ and 
$tr^q_{X\otimes Y}\,=\,tr^q_X\cdot tr^q_Y\,$, where the multiplication is in the
sense of Theorem~\ref{thm-coend=alg}.

\prg{2.3}{ Coends for Special Categories\ }
 In this paragraph we wish to give some
more reality to the formalism described above, by computing the coend explicitly
for semisimple categories and representation categories. For $\C=G-mod\,$, where
$G$ is a finite group, $F$ will be given by the space of functions on $G$ with
the coadjoint action on it. It is well known that the matrix elements of
irreducibles are a basis of this space so that by 
$$
C(G)\,=\;\bigoplus_{j\in {\cal J}_0}End_{\bf C}\bigl(V_j\bigr)^*\;\cong\;
\bigoplus_{j\in {\cal J}_0}V_j^*\otimes V_j\qquad,
$$
we can present $F$ also as a sum (or tensors) of irreducibles. 
The other cases are given by the obvious generalizations
of these two pictures.

 The answer for semisimple situation is given by the next easy result:
\blm\lll{lm-coend-ss}
Suppose \C\ is an abelian, semisimple, rigid tensor category, for which the set 
$\cal J$ of 
isomorphism classes of irreducible objects is finite. Further, let ${\cal J}_0\,$ be any set of
representatives of $\cal J$.

Then the coend from (\ref{eq-coend-not}) exists and  is given by 
$$
F\;=\;\bigoplus_{j\in{\cal J}_0}j^{\vee}\otimes j\qquad.
$$
The summand $i_X^j\,:\,X^{\vee}\otimes X\,\to\,j^{\vee}\otimes j$ of the transformation
is the canonical (identity) element in
$$
i_X^j\;\in\;Hom_{\cal C}(X^{\vee},j^{\vee})\otimes Hom_{\cal C}(X,j)\,\cong\,
End_{\bf C}\bigl(Hom_{\cal C}(j,X)\bigr)\qquad.
$$

\elm

{\em Proof: } In the last identity we used the isomorphism $Hom_{\cal C}(X^{\vee},j^{\vee})\cong Hom_{\cal C}(j,X)$ obtained from rigidity. Moreover, we used
the duality $Hom_{\cal C}(X,j)\cong Hom_{\cal C}(j,X)^*\,$ with respect to
the pairing 
$f,g\mapsto f\circ g\in End_{\cal C}(j)={\bf C}\,$, which is non-degenerate
by our definition of semisimplicity. Let us  introduce dual bases,
$\{e^{\alpha}_{jX}\}$ of $Hom_{\cal C}(j,X)\,$, and 
$\{f^{\alpha}_{jX}\}$ of $Hom_{\cal C}(X,j)\,$. It  follows from
the presentation in (\ref{eq-hom-iso}) that a natural transformation,
$\zeta\,:\,id\,\dot{\to}\,{\cal T}\,$, from the identity functor to another
functor, ${\cal T}\,:\,\C\,\to\,\C\,$, is determined entirely by the special
 morphisms $\zeta_j\,:\,j\,\to\,{\cal T}(j)\,$ on a representing set of 
irreducibles. In fact, any such set of morphisms yields by 
\beq\lll{eq-nat-ss}
\zeta_X\;=\;\sum_{j\in{\cal J}_0,\alpha}{\cal T}(e^{\alpha}_{jX})\circ\zeta_j
\circ f^{\alpha}_{jX}
\eeq
a natural transformation. Thus in the semisimple setting a morphism in 
$dn({\sf S})\,$, in the picture of coactions, is given by a set of morphisms 
$\delta_j\,:\,j\,\to\,j\otimes Z$ for $j\in{\cal J}_0\,$. 

In the case of
$F$ as above we can rewrite the transformations $i_X\,$ into coactions
$\theta_X^k\,:\,X\to X\otimes(k^{\vee}\otimes k)\,$ for 
the individual summands. They are determined by  $\theta_j^j\,=
\,coev_j\otimes \id_j\,:\,j\to j\otimes j^{\vee}\otimes j\,$ if $j\,=\,k\,$, and  
$\theta_j^k=0\,$ if $j\,\neq\,k\,$.

By rigidity we can rewrite the basic morphisms of a general coaction always 
in the form $\delta_j\,=\,
\bigl(\id_j\otimes\hat{\delta_j}\bigr)\circ\theta_j^j\,$, where 
$\hat{\delta_j}\,:\, j^{\vee}\otimes j\to Z\,$ is unique. Now, it follows
directly from (\ref{eq-nat-ss}) that 
$\delta_X\,=\,\bigl(\id_X\otimes\hat{\delta}\bigr)\circ\theta_X\,$, where 
$\hat{\delta}\,:=\,\oplus_j\hat{\delta_j}\,:\,F\,\to\,Z\,$. 
In particular, we can think of $r^Z\,=\,\hat\delta\,$ as a morphism in 
$dn({\sf S})\,$. 

Notice also that $i^j_j=\id_{j^{\vee}\otimes j}\,$ and $i^k_j=0$ if 
$k\,\neq\,j\,$ so that $r^Z\circ i_j\,=\,0\,$ for all $j\in {\cal J}_0\,$
implies that $r^Z\,=\,0\,$. Hence a morphism starting at $\lz F,i\rz\,$
is also unique in  $dn({\sf S})\,$, which completes the proof.

\hfill$\Box$
\medskip

The condition that $|\cal J|<\infty$ cannot be dropped here. In fact, we see 
from the proof that a coend
does not exist in a category, which is both  N\"otherian, i.e., the 
objects have only  finite decompositions, 
and which has infinite set $\cal J$. An example are the finite dimensional 
representations of a compact, continuous Lie-group, $G$. Here, the algebra of 
functions  on $G$ does not belong to $G-mod\,$.
Summarily, we can say that in the semisimple setting 
the condition $|{\cal J}|<\infty\,$ that entered
the construction of $\tau_{RT}\,$ in Paragraph 1.1 is the exact equivalent
of the condition on ``enough limits'' alluded to in the construction
of $\tau_L\,$ in  Paragraph 1.3. 
\bigskip

Next, let us discuss the coend for the case of a representation category.
Here, finite dimensionality of \A\  substitutes the condition 
$|{\cal J}|<\infty\,$ as a prerequisite for existence.

\blm\lll{lm-coend-tann}

Suppose that $\C=\Am\,$, where \A\ is a Hopf algebra with $dim(\A)<\infty\,$.
Then the coend from (\ref{eq-coend-not}) exists and can be chosen as
the dual linear space $\A^*\,$ endowed with the coadjoint action  as follows
$$
a.\rho\;=\;a''\lcoa\rho\rcoa S(a')\quad.
$$
Here $a\in\A\,$, $\rho\in\A^*\,$, the arrows $\lcoa$ and $\rcoa$ indicate
the canonical left and right regular actions, and we use the shorthand
$\Delta(a)=a'\otimes a''\,$ for the coproduct.

The transformations $i_X\,:\,X^{\vee}\otimes X\to F\,$ are defined from the
canonical maps  on the corresponding vector spaces. I.e., the element
 $l\otimes v$, where $v$ is a vector and $l$ a linear form, is mapped to the 
matrix element $\lz l,...v\rz\,\in\,\A^*\,$.
\elm

{\em Proof: } Let $Q$ and $Q^{\vee}$ be the projective representations as in
the end of Paragraph 2.1. For an object, $\lz Z,\xi\rz\,$, we consider the 
commutative diagram (\ref{eq-dinat-diag}) in the special case, where
$X=Y=Q\,$, and $f=\phi_a\,:\,Q\to Q\,$ is given for arbitrary $a\in \A$
by right multiplication $\phi_a(y)=ya\,$, which obviously intertwines the 
action of \A. We obtain the relation 
$$
\xi_Q(\rho\otimes a)\,=\,\xi_Q(a\lcoa\rho\otimes 1)\,=\,r^Z\bigl(i_Q(\rho\otimes a)\bigr)\qquad.
$$
Here we have set $r^Z\,:\,F\,\to\,Z\,:\;\psi\,\mapsto\xi_Q(\psi\otimes 1)\,$,
which also intertwines the action of \A. Since $i_Q$ is onto it follows
immediately that if $r^Z$ is actually a morphism in $dn({\sf S})\,$ it is 
unique.

It therefore remains to show that $r^Z\circ i_Y\,=\,\xi_Y\,$ also holds for
a general \A-module, $Y\,$. We look again at (\ref{eq-dinat-diag}), now with
$X=Q$ and $f_v\,:\,Q\,\to\,Y\,:\;a\,\to\,a.v\,$. If we check commutativity
for the special vector $l\otimes 1\in Y^{\vee}\otimes Q\,$ the assertion 
follows. \hfill$\Box$
\medskip

Another choice for the functor in (\ref{eq-S-func}) is given by 
${\sf S}(Y,X)\,=\,Y^{\vee}\odot X\,$, if we redefine ${\cal B}=\C\odot\C\,$.
Let us specify here also the coend, which  is formally given by $\F$ from 
(\ref{eq-FF}), in both considered cases:

The irreducible
objects of $\C\odot\C$ for a semisimple category are given by $i\odot j\,$,
where $i$ and $j$ run over the irreducibles of \C. The coend is thus 
given by the obvious  sum
$$
\F\,=\;\sum_{j\in{\cal J}_0}j^{\vee}\odot j\qquad ,
$$
and the dinatural transformations are completely analogous to those of 
Lemma~\ref{lm-coend-ss}. In fact, the corresponding proof is much easier  for
$\F$, since $Hom(a\odot b, c\odot d)=Hom(a,c)\otimes Hom(b,d)$,
 and implies by (\ref{eq-F=tFF}) the one for $F\,$.
\smallskip

In the Tannakian case we have $\C\odot\C\,=\,\bigl(\A\otimes\A\bigr)-mod$ ,
see [D]. As a vector space $\F$ is again $\A^*$ but now with 
$\A\otimes\A$-action given by
$$
(a\otimes b).\rho\;:=\;b\lcoa\rho\rcoa S(a)\qquad.
$$
As for the semisimple case we can use also here (\ref{eq-F=tFF}) in order to 
find an easier, though more abstract, proof of Lemma~\ref{lm-coend-tann}.
\medskip

The invariances $tr^q_X\,$ from (\ref{eq-q-tr}) are also easily identified.
For the semisimple case it is enough to consider irreducible $X$ because
$tr^q$ respects direct sums. Since the $i_j$ are simply injections of summands
the $tr^q_j\,$ are thus directly identified with the $\widetilde{coev}_j\,$'s.

For the Tannakian case $tr_X^q\,:\,1\to F\,$ can be thought of as a state on \A.
It is given by the following ``quantum-trace'':
\beq\lll{eq-trq-tann}
tr^q_X\;=\;tr^{can}_{{\sf V}_X}\rcoa G\qquad.
\eeq
Here $tr^{can}$ is the canonical trace over the vector space ${\sf V}_X$ of the
module $X$. Moreover, $G$ is the special, group like element that defines 
a balancing of \A. Equation (\ref{eq-trq-tann}) follows directly from
the relation 
\beq\lll{eq-G-coev}
\widetilde{coev}_{{\sf V}_X}\,:=\,(G\otimes \id)\circ T \circ coev_{{\sf V}_X}\,\quad,
\eeq
where $T$ is the transposition and $coev$ is the canonical coevaluation.
This relation and  its analogue for $ev_{{\sf V}_X}\,$ 
  will be important  in the combinatorial description of rigidity in  Paragraph 3.4.

\prg{2.4}{ Special Liftings, Pairings, and Modularity \ }
It will be useful and instructive to consider  a few special
situations, to which we can apply the  lifting of
 a dinatural transformation to a coend.
For the case with $Z=1$, we find the following:

\blm\lll{lm-coend-1}
For a category \C\ as above there is a canonical isomorphism,
$$
Hom_{\cal C}(F,1)\,\cong\, Nat(id_{\cal C},id_{\cal C})\qquad,
$$
 which respects the
algebra structure on each space. 

Moreover, for semisimple categories this space is one-to-one 
with the space of functions on ${\cal J}_0\,$; and for $\C=\Am\,$
it is canonically isomorphic to the center ${\tt Z}_{\cal A}\,=\,\A'\cap\A$.
\elm
 
{\em Proof: } 
An element in $Nat(id_{\cal C},id_{\cal C})\,$, i.e., 
a natural transformation of the identity 
functor, is given by system of morphisms $\delta_X\,:\,X\to X\,$ that 
commute with morphisms in \C. Hence they form an element in $dn({\sf S})\,$ in
the coaction  picture. The morphisms of the corresponding element 
$\lz 1,\xi\rz\,$ in the picture of dinatural transformations are
$\xi_X\,=\,(\id_{X^{\vee}}\otimes \delta_X)\circ ev_X\,$.
As remarked in the proof of Lemma~\ref{lm-coend-ss}, a natural transformation
is determined by the special morphisms $\delta_j\in End_{\cal C}(j)={\bf C}\,$,
for $j\in{\cal J}_0$, 
which proves the second assertion in the Lemma~\ref{lm-coend-1}.
Finally, for \C=\Am\ a linear map $l\,:F=\A^*\,\to\,{\bf C}$ is clearly
given by the evaluation on a unique element, $a_l\in\A\,$. In order for 
$l$ to be a morphism, $a_l$ must be invariant under the adjoint action of
\A. It is an elementary fact for Hopf algebras that this is equivalent to
saying that $a_l$ lies in the center.
\hfill$\Box$
\medskip

This allows us also to relate the two notions of a trace from (\ref{eq-tr})
and (\ref{eq-q-tr}) for a natural transformation, $\delta\,$, by the identity:
\beq\lll{eq-tr-qtr}
\delta\circ tr^q_X\;=\;tr_X\bigl(\delta_X\bigr)\quad .
\eeq

Let us denote  by $1^*\,:\,F\to 1\,$ the morphisms associated to the 
identity transformation of $id_{\cal C}\,$ or the dinatural transformation 
given by the evaluations $ev_X\,$. It corresponds to the unit element in \A,
and to the constant function with value 1 for semisimple categories.
\medskip

Next, suppose that for \C\ as above we have an element in
$dn({\sf S})\,$ given by transformations 
$\xi_X\,:\,X^{\vee}\otimes X\,\to\,A\otimes B\otimes C^{\vee}\,$. 
Using rigidity we can rewrite them as morphisms of the form
\beq\lll{eq-gen-dinat}
\hat\xi_X\,:\;A^{\vee}\otimes X^{\vee}\otimes X\otimes C\,\longrightarrow\,B\quad.
\eeq
The dinaturality condition from (\ref{eq-dinat-diag}) can thus  be reformulated 
more generally for the transformations $\hat\xi\,$. Conversely, any set of 
morphisms of this form that satisfy this general dinaturality can be lifted to
the coend. I.e., we find a morphism 
$r^{ABC}\,:\,A^{\vee}\otimes F\otimes C\,\to\,B\,$, such that we have the
factorization $\hat\xi_X\,=\,r^{ABC}\circ(\id_{A^{\vee}}\otimes i_X\otimes 
\id_C)\,$. For balanced or braided categories we also do not have to assume 
that the first tensor factor in (\ref{eq-gen-dinat}) is in fact a conjugate. 
\medskip

The universality of the coend  is also 
used in [L] to define the morphisms of a categorical Hopf algebra. For example 
the morphisms 
\beq\lll{eq-def-coprod}
\Delta_X\,:\;X^{\vee}\otimes X\,\TO{$\id_{X^{\vee}}\otimes coev_X\otimes \id_X$}{50}\,
X^{\vee}\otimes X\otimes X^{\vee}\otimes X\,\TO{$i_X\otimes i_X$}{35}\,F\otimes F
\eeq
are dinatural and thus lift to a coproduct $\Delta\,:\,F\to F\otimes F\,$.
\medskip 

As another application  let us explain the condition of modularity (in its strongest
version)  as the non-degeneracy of a pairing, $\omega\,$, of the coend.
It is constructed from the monodromy transformation 
$\gamma(X,Y):=\epsilon (Y,X)\epsilon (X,Y)\,\in\,End_{\cal C}(X\otimes Y)\,$,
where $\epsilon (X,Y)\,:\,X\otimes Y\to Y\otimes X\,$ is the braid constraint
of \C. Specifically, we consider the morphisms
\beq\lll{eq-omega-XY}
\omega_{X,Y}\,:\;X^{\vee}\otimes X\otimes Y^{\vee}\otimes Y\,\stackrel 
{\,\id_{X^{\vee}}\otimes \gamma(X,Y^{\vee})\otimes\id_Y\,}{\hbox to 40pt{\rightarrowfill}}\,X^{\vee}\otimes X\otimes Y^{\vee}\otimes Y\,
\stackrel {\,ev_X\otimes ev_Y\,}{\hbox to 30pt{\rightarrowfill}}\;1\quad.
\eeq

For fixed $Y$ this is  a system of morphisms as in 
(\ref{eq-gen-dinat}) with  $A,B=1\,$, and $C=Y^{\vee}\otimes Y\,$. 
The corresponding dinaturality condition follows from the naturality of 
$\gamma$ so that we can lift it to a map 
$\omega_Y\,:\,F\otimes Y^{\vee}\otimes Y\,\to\,1\,$. Now, for every $X$
the $\omega_{X,Y}$ are dinatural also with respect to $Y$. From the 
uniqueness of liftings it thus follows that $\omega_Y$ shares this property,
and therefore lifts itself to a morphism 
\beq\lll{eq-omega-coend}
\omega\,:\,F\otimes F\,\to\,1\qquad.
\eeq
Summarily, we find that $\omega_{X,Y}\,=\,\omega\circ(i_X\otimes i_Y)\,$.
An important property that was already implied in Theorem~\ref{thm-coend=alg}
is given in the next lemma.
\blm\lll{eq-pair-Hopf}
The pairing $\omega$ from (\ref{eq-omega-coend}) 
 is a {\em Hopf pairing}, i.e., $m:F\otimes F\to F$ and
$\Delta:F \to F\otimes F\,$ are adjoints of each other.
\elm

We now define a category with a coend to be {\em modular}
 if the pairing $\omega$
is non-degenerate. To illustrate the meaning of this condition let us
remark that a symmetric category is never modular unless it is $\cong\vc\,$,
for which $F={\bf C}\,$.
This follows since $\gamma(X,Y)=1$ and hence $\omega=1^*\otimes 1^*\,$, 
which is of rank one. 
In this sense modularity expresses the fact that all objects have 
non-trivial braiding.

For \C=\Am\ it can be easily worked out that the corresponding linear map 
$\A^*\otimes\A^*\to{\bf C}\,$ is given by evaluation on the element
\beq\lll{eq-omega-A}
\omega\;=\;1\otimes S({\cal R}'{\cal R})\;=\;\sum_{ij}f_ie_j
\otimes S(e_if_j)\quad.
\eeq
Here ${\cal R}=\sum_{j}e_j\otimes f_j\,\in\,\A\otimes\A\,$ is the so called
$R$-matrix of the quantum-group \A. An interesting special case is the double
$\A\,=\,D({\cal B})\,$ of a Hopf algebra $\cal B$. It is characterized by the properties
that \A contains both 
$\cal B$ and $\Bo$ (where the latter is ${\cal B}^*\,$ with opposite comultiplication)
as sub Hopf algebras, and that both maps
${\cal B}\otimes\Bo\,\to\,\A\,$ and 
$\Bo\otimes{\cal B}\,\to\,\A\,$, given by multiplication in \A, 
 are isomorphisms
of vector spaces. For the $R$-matrix the  $e_j$ form  a basis of 
${\cal B}\,$ and  the $f_j\in{\cal B}^*\,$ are given by the dual basis.
It is now clear from (\ref{eq-omega-A}) that doubles always give rise
to modular categories.
\medskip

In the semisimple case the pairing is given by its summands $\omega_{i,j}\,$,
i.e., the morphisms $\gamma(i,j^{\vee})\,$. 
A popular modularity condition that 
was  used, e.g., in the works of [RT],  is to require non-degeneracy
for the restriction of the
pairing to  invariance, i.e., for 
\beq\lll{eq-omega-inv}
Inv(\omega)\,:\;Inv(F)\otimes Inv(F)\to {\bf C}\;:\;f\otimes g\mapsto \omega(f\otimes g) 
\eeq
In the basis $\{\widetilde{coev}_j\}\,$ of $Inv(F)\,$ the matrix of this 
pairing is what is commonly called the $S$-matrix:
\beq\lll{eq-S-mat}
S_{i,j}\,=\,tr_{i\otimes j^{\vee}}\bigl(\gamma(i,j^{\vee})\bigr)\quad.
\eeq
The non-degeneracy of $Inv(\omega)\,$ is then reexpressed  as the 
invertibility of 
the $S$-matrix. This does in fact also imply non-degeneracy of $\omega$, but
for non-semisimple categories $Inv(\omega)$ may be degenerate although the 
category is modular, see the example at the end of Paragraph 3.5. 
In the next  paragraph
we shall find an even weaker and more natural condition 
 that is equivalent to modularity  for general categories.

\prg{2.5}{ Integrals and Cointegrals\ }
 Integrals 
play an important r\^ole in the constructions of the invariants 
$\tau_{HKR}$ and $\tau_L\,$, but they enter also cellular invariants,
as for example $\tau_{Ku}$ of
G. Kuperberg, see [Ku]. Moreover, they  are  at the center of study for
Hopf algebras themselves. 

A left integral of a Hopf algebra \A,  in
the classical sense, is defined as a linear form, $\mu\in{\cal A}^*$,
 with the property
\beq\lll{eq-int-class}
1\mu(y)\,=\,\id_{\cal A}\otimes\mu\bigl(\Delta(y)\bigr)\,
\equiv\,y'\mu(y'')\qquad{\rm for\ all\ }\;y\in\A\quad.
\eeq
It is an easy exercise to show that $\mu$ for the functions on a 
finite group $G$ is exactly the integration against the Haar-measure on $G$.
The notion of a right integral is analogous. We also have cointegrals,
which are integrals for the dual algebra $\A^*\,$. If we assume further
that they are given by an evaluation on an element $\lambda\in\A$ we have the 
following defining relation:
\beq\lll{eq-coint-class}
y\lambda\,=\,\epsilon(y)\lambda\qquad.
\eeq
Hence with $\epsilon$ being the counit, this means that 
$\lambda\in Inv_{\cal A}(Q)\,$, with $Q$ as in Paragraph 2.2.

It has been shown by  Sweedler, [Sw], they exist iff \A is finite
 dimensional, in which case they are unique up to scalars, i.e., 
$dim\bigl(Inv_{\cal A}(Q)\bigr)=1\,$.
\medskip

The notion of integrals and cointegrals can also be put in the more
general categorical language if we think of the coend as the Hopf algebra
dual to \A. Since $F$ is now an abstract object it does not make sense to
say that $\mu$ is an element of $F$, but it is still an element in 
the invariance $Inv(F)\,$. The integral and the cointegral of $F$  are thus 
defined as morphisms,
\beq\lll{eq-int-mor}
\mu\,:\;1\,\to\,F\qquad\; {\rm and}\;\qquad \lambda\,:\;F\,\to\,1\quad.
\eeq
The characterizing relations for left integrals from (\ref{eq-int-class})
and (\ref{eq-coint-class}) can now be rewritten as the following
pair of commutative diagrams, that make sense also in non-Tannakian  categories:
\beq\lll{eq-int-diag}
\bar{ccccccc}
F\,&\stackrel{\mbox{$\id_F\otimes\mu$}} 
{\hbox to 35pt{\rightarrowfill}}&\,F\otimes F\,&
\qquad\qquad& F\,&\stackrel{\mbox{$\Delta$}} 
{\hbox to 35pt{\rightarrowfill}}&\!\!\!F\otimes F\\
1^*\Biggl\downarrow\quad&&\quad\Biggl\downarrow m&&
\lambda\Biggl\downarrow\quad&&\qquad\;\Biggl\downarrow \id_F\otimes\lambda\\
1\,&\stackrel{\mbox{$\mu$}} {\hbox to 35pt{\rightarrowfill}}&F&&
1\,&\stackrel{\mbox{$i_{1}$}} {\hbox to 35pt{\rightarrowfill}}&\!\!\!F
\ear
\eeq

The existence and uniqueness of integrals, proven in the classical case in [Sw],
 has been generalized to the category situation in [L].
For general, braided  Hopf algebras $H$ in \C\ we actually only have 
$\mu:\,\alpha\to H\,$, where $\alpha$ is  an invertible object, i.e., 
$\alpha^{\vee}\otimes\alpha =1\,$. However, for the algebra obtained from 
the coend $F$ of a modular category
we obtain $\alpha=1$, and, further, that the integrals are {\em two-sided} (in the
braided setting).

In [LS] it is shown that integrals and cointegrals of a Hopf algebra
give rise to 
non-degenerate pairings. They are given by the morphisms
\beq\lll{eq-beta-pair}
\beta\,:\,F\otimes F\,\TO{$m$}{30}\,F\,\TO{$\lambda$}{30}\,1\;,\qquad\;{\rm and}\;\qquad
\beta^{\dagger}\,:\,1\,\TO{$\mu$}{30}\,F\,\TO{$\Delta$}{30}\,F\otimes F\qquad.
\eeq

In fact, $\beta$ and $\beta^{\dagger}\,$ are almost inverses of each other.
The following is a categorization of results in [LS] and [Rd]. 

\blm\lll{lm-beta-inv}
Suppose $\lambda$ and $\mu$ are a cointegral and an integral of a braided 
Hopf algebra in \C.
Then they can be normalized such that $\lambda\circ\mu\;=\;1$. 

Suppose $\lambda$ and $\mu$ are a right integral and a left cointegral in this normalization.
Then we have
$$
\Gamma^{-1}\;=\;(\id_F\otimes \beta)\circ(\beta^{\dagger}\otimes \id_F)\quad,
$$
where $\Gamma$ is the braided antipode, and the pairings are as in (\ref{eq-beta-pair}).
\elm

{\em Proof: } The fact that the composition is non-zero and thus can be normalized to 
1 follows from the Fundamental Theorem for Hopf modules, which is also used to
prove the existence of integrals, and is proven to hold for braided categories in [L].

The axioms of a braided Hopf algebra with braided antipode $\Gamma$ imply that the 
following diagram commutes:
\beq\lll{eq-pf-diag}
\bar{ccccc}
F\otimes F&\TO{$\Delta\otimes \id_F$}{50}&F\otimes F\otimes F&\TO
{$\id_F\otimes\epsilon(F,F)$}{60}&F\otimes F\otimes F\\
\id_F\otimes\Delta\Biggl\downarrow\qquad\;&&&&\qquad\quad\Biggl\downarrow m\otimes\id_F\\
F\otimes F\otimes F&\TO{$(\Delta\circ m)\otimes \Gamma$}{50}&F\otimes F\otimes F&
\TO{$\id_F\otimes m$}{60}&F\otimes F
\ear
\eeq

If we multiply the morphisms in $End_{\cal C}(F\otimes F)\,$ identified in 
(\ref{eq-pf-diag}) with $\id_F\otimes\mu$ from the left and with $\lambda\otimes\id_F$ from
the right we obtain the identity from Lemma~\ref{lm-beta-inv}.\hfill$\Box$
\medskip

Since we defined another canonical pairing in the previous paragraph we obtain
a canonical endomorphism:
\beq\lll{eq-Smat-gen}
S_F\,:\;F\,\TO{$\beta^{\dagger}\otimes \id_F$}{45}\,F\otimes F\otimes F\,\TO{$\id_F\otimes \omega$}{45}\,F\quad.
\eeq
The pairing from (\ref{eq-omega-coend}) can be written as $\omega\,=\,
\beta\circ\bigl(\id_F\otimes(\Gamma\circ S_F)\bigr)\,$, and modularity is the same
as invertibility of $S_F\,$. Also, the matrix in (\ref{eq-S-mat}) is equivalent to
 $Inv(S_F)$ in the semisimple case.

The following reformulation of modularity is motivated by topological considerations
that we shall discuss in Paragraph 3.1.

\btm\lll{thm-int-modular}
Suppose $\omega$ is a Hopf pairing of a braided, categorical Hopf algebra, $F$,
 with cointegral, $\lambda\,:\,F\to 1\,$. Then $\omega$ is non-degenerate if and only 
$\lambda$ lies in the image of $\omega$. I.e., there is $\mu\in Inv(F)$, such that 
\beq\lll{eq-lom}
\lambda\,:\;F\,\TO{$\mu\otimes\id_F$}{40}\,F\otimes F\,\TO{$\omega$}{30}\,1\qquad.
\eeq
If it exists, $\mu$ has to be an integral of $F$.
\etm

{\em Proof :} Since $\omega$ is a Hopf pairing it is clear that, if it is also
non-degenerate, the unique solution to (\ref{eq-lom}) is given by 
an integral. 
Conversely, assume we have a solution. The equation that characterizes $\omega$ as a 
Hopf pairing is
\beq\lll{eq-Hopf-pair}
\bar{ccc}
F\otimes F\otimes F\,&\TO{$\Delta\otimes\id_{F\otimes F}$}{50}\;
F\otimes F\otimes F\otimes F\;\TO{$\id_F\otimes\omega\otimes\id_F$}{50}\,&\,
F\otimes F\\
\id_F\otimes m\Biggl\downarrow\qquad\quad&&\quad\Biggl\downarrow\omega\\
F\otimes F\,&\TO{$\omega$}{200}\,&\,1
\ear\hfill .
\eeq
If we multiply to this $\mu\otimes \id_{F\otimes F}\,$ from the left, where $\mu$
is a solution to (\ref{eq-lom}), and use the definition in (\ref{eq-Smat-gen})
we arrive at
$$
\beta\;=\;\omega\circ(S_F\otimes\id_F)\qquad.
$$
This shows by Lemma~\ref{lm-beta-inv} that $\omega$ is invetible.\hfill$\Box$
\medskip

We conclude this chapter with a characterization of the integrals in the two
prominent, special cases:

By Lemma~\ref{lm-coend-1} there is a natural transformation of the identity
associated to the integral,
which is given by morphisms $\lambda_X\,:\,X\to X\,$, for every $X$. In the semisimple 
setting it is thus given by a function on ${\cal J}_0\,$. From the definition of the
coproduct in (\ref{eq-def-coprod}) and axiom (\ref{eq-int-diag}) we find that the 
numbers have to solve the equation
$\lambda_j i_j\,=\,\lambda_j (i_1\circ ev_j)\,$. Thus we have
$\lambda_j=0$ if $j\neq 1$ so that
$\lambda_X$ is by formula (\ref{eq-nat-ss}) a multiple of
the full projection onto the invariance of
$X$.

The substitute for the modularity condition from Theorem~\ref{thm-int-modular} was already 
proposed in [Ke5] in the semisimple case. It was expressed by the requirement that
$1\in im(S)\,$, where $S$ is the matrix from (\ref{eq-S-mat}), and $1$ is the vector with
$1_j\,=\,\delta_{1,j}\,$, i.e., it is proportional to $\lambda\,$. In the modular
case the preimage is necessarily given by the quantum-dimensions, i.e.,
\beq\lll{eq-Sd=1}
\sum_{j\in{\cal J}_0}S_{i,j^{\vee}}d(j)\;=\;{\cal D}^2\delta_{j,1}
\qquad{\rm with}\quad{\cal D}^2=\sum_{j\in{\cal J}_0}d(j)^2\quad.
\eeq
This is shown in  [T], but see also [Ke5] for a short computation. The integrals, related 
as in Theorem~\ref{thm-int-modular}, and with normalization $\lambda\circ\mu=1$, are thus
given by
\beq\lll{eq-int-ss}
\mu\,=\;{\cal D}^{-1}\sum_{j\in{\cal J}_0}\,d(j)\,\widetilde{coev}_j\;\equiv\;
{\cal D}^{-1}\sum_{j\in{\cal J}_0}\,d(j)\,tr^q_j\qquad\quad{\rm and}\quad\qquad 
\lambda\;=\;{\cal D}\,ev_1\;\equiv\;{\cal D}\,1\;,
\eeq
for a semisimple category with $|{\cal J}_0|<\infty\,$.
The alternate notation for the basis vectors in $Inv(F)\,$ 
 is obtained from (\ref{eq-q-tr}).
\medskip

Finally, for \C=\Am, we know that for a modular Hopf algebra the (classical) 
{\em right} integral $\mu^R\in \A^*\,$ is invariant under the adjoint action
- i.e., $\mu^R$ is  an element in 
$Hom_{\cal A}(F,1)\,$ - 
which, by relations from [Rd], is equivalent to saying that the {\em comodulus}
of \A\ is trivial. In fact, it is noticed in [FLM] and also [Ke2] that $\mu^R$
does coincide with the categorical integral of $\C=\Am$ we
have discussed above if \A\ is a double. 

Also, triviality of the comodulus  means exactly that the 
left cointegral $\lambda\in \A$ is two-sided, which in turn is the same as
saying that $\lambda$ lies in the center of \A. By Lemma~\ref{lm-coend-1} the
evaluation on $\lambda$ gives us therefore also the categorical cointegral, and 
the associated natural transformations $\lambda_X$ are given by acting with
$\lambda\in\A$ on the respective module.  For details, see, e.g., [Ke2].
\bigskip


\setcounter{chapter}{3}

\subsection*{ 3) Lyubashenko's Invariant, and Derivations of
 $\tau_{HKR}$ and $\tau_{RT}\,$\ :}\lll{pg-3}
\

In this chapter we will discuss the general construction of
three-manifold invariants starting from a modular, abelian,
braided tensor category, \C. We shall apply the special results on coends for
semisimple and Tannakian categories from the previous chapter,
in order to show that $\tau_L\,$  
specializes to the invariants of Reshetikhin-Turaev and Hennings,
respectively, thus proving   Theorem~\ref{thm-main-rel}.

\prg{3.1}{ Construction of $\tau_L\,$ }
 As for $\tau_{RT}\,$ the construction 
of $\tau_L\,$ starts by presenting a three-manifold $M$ by surgery along
a framed link $\li\subset S^3$, or, equivalently, a link of ribbons. 

The link can be thought of as a special element
in the category of isotopy classes of oriented, framed tangles, $\tgl\,$. An object is 
an ordered set of labels $\{a_1,\ldots,a_N\}\,$, and a morphism 
$\{a_1,\ldots,a_N\}\,\to\,\{b_1,\dots,b_M\}\,$ is the generic 
projection of a framed tangle in  the strip  ${\bf R}\times [0,1]\,$.
It consists of ribbons that end in $N$ intervals along the top-line 
${\bf R}\times\{1\}\,$ and in $M$ intervals along the bottom-line.
We require that the labels of the intervals, that are connected by a ribbon, 
are the same. We also admit  closed ribbons so that $\li\in  
End_{\tgl}(\emptyset)\,$ for the projection of a closed link.
\medskip

For a given abelian category \C\ we call a {\em coloring}, {\tt Col}, an
assignment, $a\mapsto X_a\,$, of labels to general objects of \C, 
together with  a sense of direction for each component of a given tangle.
To every coloring we then associate a functor:
\beq\lll{eq-fc-li-ab}
I^{\tt Col}\,:\;\tgl\,\longrightarrow \,\C\qquad.
\eeq

The construction of $I^{\tt Col}$ proceeds as in [RT]. To a label set,
$\{a_1,\ldots,a_N\}\,$, we associate the object 
$X^{\#}_{a_1}\otimes\ldots \otimes X^{\#}_{a_N}\,$, where $X^{\#}_a=X_a\,$ if the 
attached strand is directed downward and $X^{\#}_a=X^{\vee}_a\,$ if the strand
is directed upward. A tangle projection  is always isotopic to one with a generic 
height coordinate, in which the tangle is presented as the
 composite of elementary tangles. The latter are allowed to contain only one crossing, one 
maximum, or one minimum. They are mapped by $I^{\tt Col}$ to the braid 
constraint $\epsilon(X,Y)\,$ and the rigidity morphisms, respectively.
In our convention the latter assignment is as follows:
 
\beq\lll{fig-rigid}
\begin{picture}(430,70)
\put(20,40){\oval(40,40)[t]}
\put(17,57){$<$}
\put(45,50){$\mapsto coev_X$}
\put(130,60){\oval(40,40)[b]}
\put(128,37){$<$}
\put(157,50){$\mapsto ev_X$}
\put(245,40){\oval(40,40)[t]}
\put(243,57){$>$}
\put(270,50){$\mapsto \widetilde{coev}_X$}
\put(355,60){\oval(40,40)[b]}
\put(352,37){$>$}
\put(380,50){$\mapsto \widetilde{ev}_X$}
\end{picture}
\eeq

Moreover, we assign a balancing element, $v_{X^{\#}}\,$, with 
$ v\in Nat_{\cal C}(id,id)\,$   to every local $2\pi$-twist of a ribbon with 
coloring $X$.  By the framed version of the first Reidemeister move, it
relates the rigidity morphisms to their flipped counterparts.
\medskip

With $I^{\tt Col}(\emptyset)=1\,$, this procedure assigns to every link, $\li\,$, with
$N$ components and every 
coloring with objects $(X_1,\ldots,X_N)\,$ a number 
$I^{(X_1,\ldots,X_N)}\bigl(\li,\C\bigr)\,$. We could now, as for $\tau_{RT}\,$,
try to find a combination of such invariants that is invariant under 
2-handle slides. The following lemma shows that it is not possible to find 
anything new for non-semisimple categories in this way:

\blm\lll{lm-JH-fac}
To an object, $X\,$, of a N\"otherian, abelian, braided, tensor category, \C, we
associate the set ${\cal J}(X)\,$, which consists of all irreducible factors of
a Jordan-H\"older series of $X$ that have non-vanishing quantum-dimensions
(with repeated isomorphic objects). 

For any coloring we then have
$$
I^{(X_1,\ldots,X_N)}\bigl(\li,\C)\;=\;
\sum_{j_1\in{\cal J}(X_1),\ldots,j_N\in{\cal J}(X_N)}I^{(j_1,\ldots,j_N)}\bigl(\li,\overline{C}^{tr}\bigr)\qquad.
$$
\elm

This follows from the fact that for each $X_{\nu}$ the invariant can be 
expressed as the canonical trace $tr_{X_{\nu}}$ over a morphism that depends on
$X_{\nu}\,$ like a natural transformation. The value thus only depends on its
image in the semisimple trace-quotient of \C\ . The object $X_{\nu}$ itself factors
into the object $\overline{X_{\nu}}\,=\,\oplus_{j_{\nu}\in{\cal J}(X_{\nu})}\,j_{\nu}\,$
in $\overline{\C}^{tr}\,$.
\medskip

In Lemma~\ref{lm-coend-1} we learned that natural transformations of the 
identity, as in the computation alluded to above, are naturally identified 
with a coinvariance-morphism  or state on the coend.  In order to obtain such
a state from a given link, $\li\,$, we have to open the components of $\li$
in a controlled way. The procedure can be described  as follows:

We first introduce  a horizontal line ${\bf R}_t$ above a region in the plane
that contains a projection of $\li\,$. For each component, 
$C\cong S^1\times [0,1]\,$, of the link we
introduce a  ribbon, $R_C$, that starts at  an interval in ${\bf R}_t$ and ends
in an interval on $\partial C$. Then we split $R_C$ down the middle so that
we obtain two parallel  ribbons, $R_C'$ and $R_C''\,$. We also cut the component
$C\,$, where $R_C$ has been attached, as indicated in the diagram~\ref{fig-rib-split}  

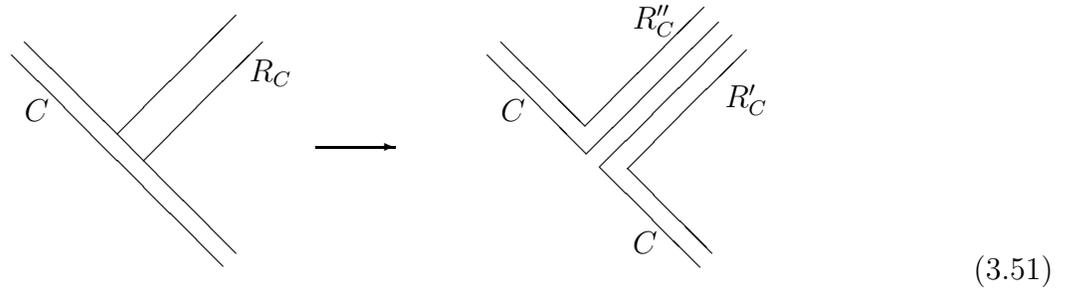
\beq\lll{fig-rib-split}
\begin{picture}(300,100)(0,0)
\put(10,90){\line(1,-1){80}}
\put(5,85){\line(1,-1){80}}
\put(45,55){\line(1,1){45}}
\put(55,45){\line(1,1){45}}
\put(10,60){$C$}
\put(95,75){$R_C$}
\put(120,50){\vector(1,0){30}}
\put(190,90){\line(1,-1){32}}
\put(185,85){\line(1,-1){37.5}}
\put(227.5,42.5){\line(1,-1){38}}
\put(238,42){\line(1,-1){32}}
\put(222,58){\line(1,1){45}}
\put(222.5,47.5){\line(1,1){50}}
\put(227.5,42.5){\line(1,1){50}}
\put(238,42){\line(1,1){45}}
\put(190,60){$C$}
\put(240,10){$C$}
\put(275,65){$R_C'$}
\put(240,95){$R_C''$}
\end{picture}
\eeq

Thus we obtain a tangle, $\li^{split}\,$, which contains no closed components
anymore, but only ribbons that start and end in ${\bf R}_t$. Moreover, the
start and end points are always neighbors. Hence, as a morphism 
$\{a_1,\ldots,a_{2N}\}\,\to\,\emptyset\,$ in $\tgl$, where $N$ is the number of components of 
$\li$, we have for the labels $a_{2i-1}=a_{2i}\,$.

As described in [KL]  this tangle represents
the manifold $M^{split}$, defined in the introduction of Chapter 1, 
(see also [Ke5] of Paragraph 3.5).

Let us choose directions such that the
strand at $a_{2i}\,$ is going downwards. A coloring, {\tt Col},
 is now uniquely given by the  $N$ objects from $a_{2i}\mapsto X_i\,$.  
Applying the functor 
$I^{\tt Col}\,$ from (\ref{eq-fc-li-ab}), we obtain morphisms
\beq\lll{eq-I-X-1}
I^{(X_1,\ldots,X_N)}\bigl(\li^{split},\C\bigr)\,:\;X_1^{\vee}\otimes X_1\otimes\ldots\otimes X^{\vee}_N\otimes X_N\;\longrightarrow\;1\qquad.
\eeq
These maps are dinatural in any of the arguments. This means for example that
for fixed $X_1,\ldots,X_{j-1},X_{j+1},\ldots,X_N\,$ 
they have the property of the morphisms  described in (\ref{eq-gen-dinat}) with
respect to $X=X_j$, setting
$A^{\vee}=X_1^{\vee}\otimes X_1\otimes \ldots\otimes X_{j-1}\,$,
$C=X_{j+1}^{\vee}\otimes \ldots\otimes X_N\,$, and $B=1\,$. This follows from the
fact that all elementary morphisms are natural so that a morphism 
$f:X_j\to Y_j\,$ can be pushed through an entire component of a diagram. The transformation
can thus be lifted to the coend. Repeated application yields a state on the $N$-fold tensor product of the coend, i.e., a  morphism 
$\Ii\bigl(\li^{split},\,\C\bigr)\,:\;F^{\otimes N}\,\to\,1\,$, from which we recover the ones in
(\ref{eq-I-X-1}) by multiplying the transformations $i_{X_j}\,$.
The question how the morphism $\Ii\,$ relates to  invariants is answered 
next:

\blm\lll{lm-int-inv}
For a sequence of invariances, $\mu_j\in Hom_{\cal C}(1,F)\,$,  define a number by
\beq\lll{eq-int-inv}
\tau_L\bigl(\li^{split},\,\C,\,\vec{\mu})\;:=\;\Ii\bigl(\li^{split},\,\C\bigr)\circ(\mu_1\otimes\ldots\otimes\mu_N)\quad.
\eeq
\ben
\item If all of the $\mu_j$ are invariant under the braided antipode $\Gamma$
then $\tau_L$ only depends on $\li$, i.e., not on the choice of splitting ribbons
$R_C\,$.
\item Suppose $\mu_j\,=\,tr^q_{X_j}$, with 
$tr^q_X\,$ as in (\ref{eq-q-tr}). Then $\tau_L$
is identical to the invariant from Lemma~\ref{lm-JH-fac}.
\item If $\mu_j=\mu$ are integrals of $F$, then $\tau_L$ is invariant
under 2-handle slides. 
\een
\elm

{\em Proof: } In order to verify independence of the choice of the splitting
ribbons 
we have to be able to switch an under-crossing with $R_C$ to an over-crossings,
which is possible since all $\mu_j$ are in the invariance. It remains to ensure
that we can switch the attachment from one side of $C$ to the other. This is
taken care of by the antipode. In a precise version of this fact we also 
have to take care of orientations of $R_C$ as a surface with boundary, see [L]
or [KL] for details.

It is easy to see that we can recover $\li$ up to isotopy from $\li^{split}\,$
by placing $N$ small, flat arcs on top of ${\bf R}_t\,$, joining an interval with label $a_{2i-1}\,$  to an interval with label $a_{2i}\,$. With given colorings
and the conventions from (\ref{fig-rigid}) this corresponds to the insertion 
of flipped coevaluations for the $\Ii$ factored into the respective 
transformation.

For the fact that the axiom on the left of (\ref{eq-int-diag}) directly implies
the 2-handle slide, see [L], but also [KR] in the Hennings-framework,
\hfill$\Box$
\medskip

The invariant from  Part $3$ of the lemma is basically Lyubashenko's 
invariant, 
$\tau_L\,$. The only thing  we have omitted from the discussion
 is the invariance under the
signature-move (or ${\cal O}_1$-move in [Ki]). 
As we remarked in Paragraph 1.4, our notion of  TQFT's is really defined 
for a central extension of the cobordisms
categories by $\Omega_4\,$.
The ${\cal O}_1$-move is exactly the generator of this
extension. We shall therefore substitute it by a move under which two link
diagrams are equivalent if they not only describe the same three-fold 
but also have linking-matrices with the same signature, $\sigma(\li)\,$.
It is given by the addition or removal of an isolated Hopf-link, $\nhopf$,
for which
one component is 0-framed, see also the $\eta$-move in [Ke5]. Due to 
multiplicativeness of $\tau_L$ and 2-handle slides over the 0-frame,
it is thus enough to check that
\beq\lll{eq-eta-inv}
\tau_L\bigl(\hopf\bigr)\;=\;\tau_L\bigl(\twhopf\bigr)\;=\;1\qquad.
\eeq
In a variation of Kirby's calculus described next, it turns out that 
this imposes only a normalization condition on $\mu$, which determines
the integral up to a sign. We will not consider the usual renormalizations 
of the entire invariant, e.g., by 
$\tau_L\bigl(\bigcirc^{\pm 1}\bigr)^{\mp\sigma({\cal L})}\,$, required by the 
original ${\cal O}_1$-moves, or by $\tau_L\bigl(S^1\times S^2\bigr)^{-1}$,
which is not defined in the non-semisimple case.
\medskip

Let us outline here the mentioned
 modification of Kirby's calculus, namely the so-called
``bridged link calculus'' from [Ke5]. It has the  advantage that it involves 
only local moves, at the price of introducing additional data, that 
indicates  1-handle surgeries. In one  version of the calculus 
this data is represented 
by horizontal coupons, which may have strands $R_1,...,R_n\,$ 
entering at the top side and leaving at the 
bottom side, and is sometimes also expressed by  a ``Kirby's dotted circle'' around 
$R_1,...,R_n\,$. The relevant moves of this calculus are depicted in
the following diagram. We assume as usual blackboard-framing, i.e.,  
flat ribbon-projections:

\beq\lll{fig-BL-mvs}
\begin{picture}(340,100)(20,0)
\put(17,90){$X$}
\put(20,44){\line(0,1){38}}
\put(20,50){\oval(32,18)[b]}
\put(18,50){\oval(28,18)[tl]}
\put(22,50){\oval(28,18)[tr]}
\put(20,38){\line(0,-1){25}}

\put(50,50){$=$}

\put(80,90){$X$}
\put(87,58){\line(0,1){24}}
\put(87,42){\line(0,-1){29}}
\put(74,42){\framebox(26,16){$\lambda_X$}}

\put(20,-10){\em Modification}

\put(150,42){\framebox(26,16){}}
\put(180,58){\oval(34,40)[t]}
\put(180,42){\oval(34,40)[b]}
\put(197,42){\line(0,1){16}}
\put(205,50){$=\;\emptyset$}

\put(150,-10){\em 1-2-Cancellation}

\put(276,90){$X$}
\put(283,68){\line(0,1){14}}
\put(283,52){\line(0,-1){43}}
\put(270,52){\framebox(26,16){$\lambda_X$}}
\put(285,30){\line(1,1){17}}
\put(281,26){\line(-1,-1){17}}

\put(310,50){$=$}

\put(336,90){$X$}
\put(343,68){\line(0,1){14}}
\put(343,52){\line(0,-1){22}}
\put(343,26){\line(0,-1){17}}
\put(330,52){\framebox(26,16){$\lambda_X$}}
\put(362,47){\line(-1,-1){36}}

\put(280,-10){\em Coupon-Slide}
\end{picture}
\eeq
\bigskip 

In the diagrams we already indicated an extension of the construction of
$\tau_L$ to links with coupons: If the strands going through a coupon from
top to bottom are colored with objects $X_1,\ldots,X_n\,$, we insert here a
morphism $\lambda_X\,:\,X\to X\,$ with $X=X_1\otimes\ldots\otimes X_n\,$.
In particular, the strands for the modification and the coupon-slide
may be replaced by  several parallel strands. The modification-move also
imposes that the $\lambda_X\,$ do belong to a natural transformation of the
identity. In addition to the moves in (\ref{fig-BL-mvs}) we also have to 
require invariance of the balancing, $v_X\lambda_X=\lambda_X\,$, in order
to guarantee invariance with respect to collective $2\pi$-twists of the 
strands entering a coupon.  It is clear that these moves also imply
the normalization conditions in (\ref{eq-eta-inv}).

In view of our discussion in Paragraph 2.5 it is not difficult to see that
a  natural candidate for $\lambda$ is  the cointegral of $F$. It is 
has become evident in Hopf algebra theory that the interplay of integrals 
and cointegrals can often be very useful in the organization and 
clarification of  involved structures and computations. It is in this sense,
the formalism  of bridged links enables us to 
apply the same principles in a purely topological setting to the relations 
between index 1 and index 2 surgeries on a three-fold.

If we are interested only in closed manifolds modularity, and the fact that
$\mu$ and $\lambda$ are integrals are not absolutely necessary. We could
try to 
define an invariant from general elements
\beq\lll{eq-adm-pair}
\mu\in Hom_{\cal C}(1,F)\qquad\qquad\quad\lambda\in Hom_{\cal C}(F,1)\;.
\eeq
Invariance under the special 1-2-cancellation from the above diagram
 imposes the constraint $\lambda\circ\mu\,=\,1\,$. Moreover, the equation 
related to the  modification-move holds for every $X$ iff the pairing 
$\omega$ maps $\mu$ to $\lambda$ as in Theorem~\ref{thm-int-modular}.
Let us call here a pair as in (\ref{eq-adm-pair}) with these two properties
{\em dual states}.

To an invariance, $\mu$, we can compute the difference $D_{\mu}\,$ of the
two morphisms from $F$ to itself in the left diagram of (\ref{eq-int-diag}).
 In the same way we can define $D^{\lambda}\,$ from the diagram for 
cointegrals. It is clear  that $\mu$ is an integral iff $D_{\mu}=0\,$,
and $\lambda\,$ analogously. In absence of strict modularity we can
also introduce the weaker notion of $\omega$-integrals and $\omega$-cointegrals
defined by the conditions
$$
\omega\circ(\id_F\otimes D_{\mu})\,=\,0\qquad{\rm and}\qquad
\omega\circ(\id_F\otimes D^{\lambda})\,=\,0\;,
$$
respectively. The first two assertions of the following lemma are
 immediately obtained from liftings to coends and and the definitions of $F$
as an algebra.
In the last one we use that $\omega$ is a Hopf pairing in a similar way as
in the proof of Theorem~\ref{thm-int-modular}.

\blm\lll{lm-omega-ints}
Suppose  $\mu$ and $\lambda$ are as in (\ref{eq-adm-pair}). 
\ben
\item The 2-handle slide of an open strand with color $X$ over a closed one 
holds for all $X$, if and only if  $\mu$ is an $\omega$-integral.

\item The coupon-slide as in diagram (\ref{fig-BL-mvs}) holds for all colorings,
if and only if $\lambda$ is an $\omega$-cointegral.

\item Suppose that $\lambda$ and $\mu$ are dual states. Then $\mu$ is an
$\omega$-integral if and only if $\lambda$ is an $\omega$-cointegral.
\een
\elm

Let us remark here that we should have anticipated the last (purely algebraic)
statement already
from the first two, using the close correspondence between algebraic and
topological structures. Specifically, this was to be expected since
introducing  the 2-handle slide in addition to the first two moves 
leads us to an equivalent calculus of surgery diagrams as introducing the
coupon-slide. Finally, it is not hard to see that, with $\lambda$ and $\mu$
as in Lemma~\ref{lm-omega-ints}, the balancing $v:F\to 1$ is in the image 
of $\omega$ so that $v\cdot\lambda\,=\,(v\circ i_1)\lambda\,=\,v_1\lambda\,$.
We thus have invariance under $2\pi$-twists at coupons if $v_1=1\,$.

In the treatment of [L] another algebra, {\bf f}, is introduced by dividing
the coend $F$ by the null-space of $\omega\,$. This allows us to construct
also ordinary TQFT's for categories with much weaker modularity conditions.
We shall, however, focus in the following on the strictly modular case,
where we can hope to represent extended structures, as discussed in Chapter 1.

\prg{3.2}{ Specialization to  $\tau_{RT}\,$ \ }
The proof of the first part of Theorem~\ref{thm-main-rel} is now a 
straightforward application of  the theory developed so far:

We insert the formula (\ref{eq-int-ss}) for the integral of the semisimple
category into (\ref{eq-int-inv}) from Lemma~\ref{lm-int-inv}), and 
reorganize the summation in the following form.
$$
\tau_L(\li,\C)\;=\;{\cal D}^{-N}
\sum_{j_1,\ldots,j_N\in{\cal J}_0}\Pi_{\nu=1}^Nd(j_{\nu})
\,\,\Ii(\li^{split},\C\bigr)\circ\bigl(tr^q_{j_1}\otimes\ldots tr^q_{j_N})\quad\;,
$$
where $\cal D$ is as in (\ref{eq-Sd=1}). Next, it follows from  Part {\em 3} of Lemma~\ref{lm-int-inv}) 
that the composite of morphisms in this sum is identical to 
$I^{(j_1,\ldots,j_N)}(\li,\C)\,$. If we set
$$
w(j_1,\ldots,j_N)\;=\;{\cal D}^{-N}\Pi_{\nu=1}^Nd(j_{\nu})\;,
$$ 
we finally obtain the $\tau_{RT}\,$-invariant in the form given in
(\ref{eq-RT-sum}). 
In this normalization $\tau_{RT}$ is strictly speaking 
an invariant of pairs, given by a three-fold together
with a signature. 
It was already observed in [RT] that there is basically no 
freedom in the choice of the weights, except for overall scalings.
>From the point of view we have taken here, this is simply a reflection of
the uniqueness of the integral of a categorical Hopf algebra.
\medskip

In the introduction we already pointed out that the $\tau_{RT}$-invariant
is often also constructed for a non-semisimple category, \C, like that of  
representations of a quantum group at a root of unity, by considering the
semisimple trace-subquotient $\overline{\C}^{tr}\,$ instead. 
Yet, by virtue of Lemma~\ref{lm-JH-fac}  there is also a way to compute 
it from $\tau_L$ for the original $\C$  with some modifications.

To be more precise let us pick a set of objects,
${\cal J}_0^*\,$, of \C, which factors into a  set of representatives of 
irreducibles, 
${\cal J}_0\,$, for $\overline{\C}^{tr}\,$. It is not hard to see from the 
definition of the functor $\C\,\to\overline{\C}^{tr}\,$ that these objects can
always be chosen as irreducibles themselves, and that any two choices of
irreducibles are isomorphic in \C. From Lemma~\ref{lm-JH-fac} we have 
in particular for $j_{\nu}\in{\cal J}_0^*\,$ that 
$I^{(j_1,\ldots,j_N)}\bigl(\li,\C\bigr)\,=\,I^{(\overline{j_1},\ldots,
\overline{j_N})}\bigl(\li, \overline{\C}^{tr}\bigr)\,$.

Let us also introduce the following invariance of the coend $F$ in \C\ :
\beq\lll{eq-mu-Css}
\mu^{semis}\;:=\;{\cal D}^{-1}\,\sum_{j\in{\cal J}_0^*\,}\,d(j)\,tr^q_j\;,
\eeq
with $tr^q$ as in (\ref{eq-q-tr}). Since $tr^q_X$ only depends on the 
equivalence class of $X$ in \C,  $\mu^{semis}\,$ is in fact a canonical
invariance, independent of the choice of ${\cal J}_0^*\,$.
\medskip

Using Part {\em 2} of Lemma~\ref{lm-int-inv}, it is now clear that we obtain
the $\tau_{RT}\,$-invariant if we substitute $\mu^{semis}\,$ for the integral
$\mu$ of \C. More precisely:
\beq\lll{eq-rrr}
\tau_{RT}\bigl(\li,\overline{\C}^{tr}\bigr)\;=\;
\tau_L(\li,\C, \mu^{semis}\bigr)\qquad.
\eeq

It is often convenient to express this modification of $\tau_L$ in
a slightly different way. To this end we remark that the 
non-degenerate invariant form defined in (\ref{eq-beta-pair}) associates
to the canonical invariance $\mu^{semis}$ a unique, canonical state,
${\sf Q}\,:\,F\to 1\,$, by the condition 
\beq\lll{def-Q-bms}
\mu^{semis}\,:\;1\,\TO{$\beta^{\dagger}$}{50}\,F\otimes F\,
\TO{$\id_F\otimes{\sf Q}$}{60}\,1\;.
\eeq
Quite obviously we have that \C\ is semisimple iff (up to rescalings)
${\sf Q}=1^*$.
Also, we find that {\sf Q} is invariant under the braided antipode $\Gamma\,$.
\smallskip

Now, in view of  the correspondence of Lemma~\ref{lm-coend-1} we can think of 
{\sf Q} also as a natural transformation of the identity, i.e., a set of
morphisms ${\sf Q}_X\,:\,X\to X\,$. 
Semisimplicity of \C\ 
is then equivalent to the condition ${\sf Q}_X\,=\,\id_X\,,\;\forall\,X\,$.
\medskip

In order to relate the substitution of (\ref{eq-rrr}) to  ``{\sf Q}-insertions'',
let us  describe next a generalization of the  procedure from Paragraph 3.1, 
that associates an  invariant 
$\tau_L\bigl(\li^{\sf Q},\C)\,$, to a link,  $\li$, with inserted 
 natural transformations. More precisely,
we consider {\em decorated} links for which a strand may  contain a chip
labeled by a natural transformation of the identity, as depicted 
in the diagram (\ref{fig-Q-ins}) below.

\beq\lll{fig-Q-ins}
\begin{picture}(40,90)(35,0)
\put(20,45){\circle{28}}
\put(20,59){\line(0,1){18}}
\put(20,31){\line(0,-1){24}}
\put(14.5,41.5){${\sf Q}_X$}
\put(16,84){$X$}
\end{picture}
\eeq

The rule to compute $I^{\tt  col}\,$ for a given coloring is now to 
insert for a chip at a downward strand with coloring $X$ the morphism
${\sf Q}_X\,$.  
To be precise we also have to require 
\beq\lll{eq-jjj}
{\sf Q}_X^t\;=\;{\sf Q}_{X^{\vee}}\qquad,
\eeq
in order for this to be well defined with respect to sliding a chip
through extrema. Note, that coupons  are also decorations
in this sense, since the cointegral in our situation is two-sided.

Let us denote by $\li^{\sf Q}\,$ the link, where we have inserted exactly 
one chip into every component of $\li$, and labeled them all with the same
{\sf Q}. In a splitting the chips can all be moved to the top-line ${\bf R}_t\,$
at the intervals with labels $a_{2i}\,$.
\medskip

We now have the following modification of (\ref{eq-rrr}):

\btm\lll{thm-RT-LQ}
Suppose \C\ is a modular category with modular quotient $\overline{\C}^{tr}$
and let {\sf Q} be  as defined in   (\ref{def-Q-bms}). 

Then
\beq\lll{eq-RT-LQ}
\tau_{RT}\bigl(\li,\overline{\C}^{tr}\bigr)\;=\;\tau_L\bigl(\li^{\sf Q},\C\bigr)\qquad,
\eeq
where the chosen integrals are those of the coends of  the respective categories.
\etm

{\em Proof:} 
Suppose $q,q'\,:\,F\to Z\,$ are related by $q\circ i_X\bigr(\id_{X^{\vee}}\otimes 
{\sf Q}_X\bigl)\,=\,q'\circ i_X\,$. Then we find from the definition of the 
coproduct in (\ref{eq-def-coprod}) that 
$q'\,=\,(q\otimes{\sf Q})\circ\Delta\,$. Hence, with 
(\ref{def-Q-bms}) we find $ q'\circ\mu\,=\, q\circ\mu^{semis}\,$. If we substitute 
for $q$ the lifted transformation of $\Ii(\li^{split})\,$ we see from the assumptions that
the lift, $\Ii\bigl((\li^{\sf Q})^{split}\bigr)\,$, where we inserted into the
transformations from (\ref{eq-I-X-1}) a ${\sf Q}_X\,$ at each second interval at
the top-line, corresponds to $q'$. A comparison with (\ref{eq-rrr}) then completes
the proof.\hfill$\Box$  
\medskip

An interesting situation would arise if we found for some category 
a (split) knot,  $\cal K$, such that  ${\sf Q}\,=\,\Ii\bigl({\cal K}^{split}, \C\bigr)\,$. In this case we can switch between the two invariants 
by connected summing of $\cal K$ to every component of $\li\,$, i.e.,
$\tau_{RT}\bigl(\li,\overline{\C}^{tr}\bigr)\,=\,
\tau_L\bigl(\li\#({\cal K}^N),\C\bigr)\,$. A more realistic question is, 
whether we can find at least a formal linear combination of knots that yields the 
canonical transformation {\sf Q}.
 
In Paragraph 3.4 we will  discuss the above  correspondence between invariants via {\sf Q}
in the more concrete 
context of Hennings-invariants.

\prg{3.3}{ Categories of Singular Tangles and Hennings' Rules}
The basic braid and rigidity morphisms of a strict, Tannakian, braided
tensor category, $\C=\Am$, can be identified with the action of
certain elements of the quantum-group \A\ on given modules. 
Specifically, these are 
an $R$-matrix and a balancing element,
\beq\lll{eq-spec-elem}
{\cal R}\,=\,\sum_je_j\otimes f_j\;\in {\cal A}\otimes{\cal A}\,,
\qquad{\rm and}\qquad G\,\in\A\qquad .
\eeq

It is a quite old idea to get rid of the modules in the computation of tangle 
invariants and try
to find a combinatorics of these elements alone. The first such attempt 
was made by Reshetikhin in  [Re], where he associated to open 
(one-component) tangles 
an element in the center ${\tt Z}_{\cal A}\,$ of \A. 
At about the same time Hennings formulated
the same  combinatorial rules, but also included a (classical) 
right integral, $\mu^R\,$,  of \A\ in the picture. This allowed him to treat
closed components of a tangle by evaluating the elements on this component
against $\mu^R$. The number that is eventually computed for a closed link 
by this procedure is the three-manifold invariant, $\tau_{HKR}\,$.
\medskip

In this paragraph we shall briefly review these rules,
first without integrals, 
in the improved version of Kauffman and Radford, using a slightly more 
categorical language. The central gadget to be defined and studied is 
a  category, \stA ,  of ``singular \A-tangles'', associated to a balanced 
Hopf-algebra \A.

\stA\ is defined quite analogously to \tgl. Its objects are again strings of labels and a morphism is again given by generic strands in a strip, that are either
closed paths, or that 
start and  end at boundary points of the strip with the same  labels.
They differ from the tangles in \tgl\ in that we do not distinguish
between under- and over-crossings,  and that we can have 
blobs sitting on a strand (similar to the chips),
which are labeled by elements in \A.
Thus the category is generated by the following, elementary morphisms:

\beq\lll{fig-sing-gen}
\begin{picture}(430,40)
\put(20,10){\oval(40,40)[t]}
\put(100,30){\oval(40,40)[b]}
\put(190,35){\line(1,-1) {30}}
\put(190,5){\line(1,1) {30}}
\put(290,5){\line(0,1){30}}
\put(290,20){\circle*{10}}
\put(297,21){$x$}
\put(335,5){for any  $x\in {\cal A}\quad$.}
\end{picture}
\eeq

Any other morphism can be obtained by composition, i.e., putting 
tangles on top of each other, by tensor products, i.e., juxtapositions of
morphisms, and, finally, by (distributive) linear combinations of 
diagrams for a fixed topological
tangle but different decorations by elements in \A. The latter shall be
indicated in a diagram by writing summation indices to the elements.
\medskip

Given  these generators we can formally characterize the category \stA\ 
by a list of relations among the generating morphisms. In \tgl\ these
had been the three Reidemeister moves (the first with framing), and,
further, the two moves due to a preferred vertical direction, namely the
extrema-cancellations and the ``crossing-symmetry''. Their singular versions
shall also apply to parts of tangles that carry no Hopf algebra elements.
Only the first  Reidemeister move shall be replaced by the move indicated
on the left of diagram (\ref{fig-st-isom}), where $G\in\A$ is a-priori
any invertible element.  An analogous relation at a maximum can 
be derived from this using the extrema-cancellation moves.  

\beq\lll{fig-st-isom}
\begin{picture}(400,100)(0,-15)
\put(30,30){\oval(40,40)[b]}
\put(10,30){\line(2,5){4}}
\put(50,30){\line(-2,5){4}}
\put(14,40){\line(1,1){36}}
\put(46,40){\line(-1,1){36}}

\put(72,43){$=$}

\put(120,30){\oval(40,40)[b]}
\put(100,30){\line(0,1){46}}
\put(140,30){\line(0,1){46}}
\put(140,55){\circle*{10}}
\put(147,56){$G$}

\put(60,-15){\em Rm1}

\put(220,40){\line(0,1){20}}
\put(250,40){\line(0,1){20}}
\put(220,60){\line(2,5){4}}
\put(250,60){\line(-2,5){4}}
\put(220,40){\line(2,-5){4}}
\put(250,40){\line(-2,-5){4}}
\put(224,70){\line(1,1){26}}
\put(246,70){\line(-1,1){26}}
\put(224,30){\line(1,-1){26}}
\put(246,30){\line(-1,-1){26}}

\put(270,50){$=$}

\put(290,4){\line(0,1){92}}
\put(320,4){\line(0,1){92}}

\put(263,-15){\em Rm2}

\end{picture} 
\eeq 
 
We also depicted here the form of the second Reidemeister move. Since we
have singular crossings it implies that we have a presentation 
of the symmetric groups instead of the braid groups as in \tgl .

Next, we give relations that  allow us to move and combine  
elements along a component of a singular tangle:

\beq\lll{eq-elem-mv}
\begin{picture}(440,85)(0,-15)
\put(20,5){\line(0,1){60}}
\put(20,20){\circle*{10}}
\put(27,21){$x$}
\put(20,50){\circle*{10}}
\put(27,51){$y$}

\put(37,35){$=$}

\put(60,5){\line(0,1){60}}

\put(60,35){\circle*{10}}
\put(67,36){$xy$}

\put(33,-15){\em Me1}

\put(135,35){\oval(40,40)[b]}
\put(115,35){\line(0,1){20}}
\put(155,35){\line(0,1){20}}
\put(155,42){\circle*{10}}
\put(162,43){$S(x)$}

\put(190,33){$=$}

\put(235,35){\oval(40,40)[b]}
\put(215,35){\line(0,1){20}}
\put(255,35){\line(0,1){20}}
\put(215,42){\circle*{10}}
\put(222,43){$x$}

\put(185,-15){\em Me2}

\put(300,55){\line(1,-1){40}}
\put(300,15){\line(1,1){40}}
\put(330,25){\circle*{10}}
\put(340,25){$x$}

\put(358,35){$=$}

\put(380,55){\line(1,-1){40}}
\put(380,15){\line(1,1){40}}
\put(390,45){\circle*{10}}
\put(397,48){$x$}

\put(352,-15){\em Me3}

\end{picture}
\eeq

Again, using extrema cancellations and crossing-symmetry, we can derive 
analogous relations  for moving an element through a maximum or
along the opposite strand of the crossing. If we have that $G$ 
implements the antipode $S$, i.e., 

\beq\lll{eq-G-S}
S^2(x)\;=\;GxG^{-1}\qquad\forall\,x\in\A\quad,
\eeq

then an element $x$ can be moved through a piece of a tangle and emerge
as $S^p(x)\,$, where the order $p$, by which the antipode acts, is invariant
under the isotopy moves. Moreover, $p$ is even if the element leaves
the piece in the same direction as it entered, and  odd if it reverses
the direction. Also, two applications of a crossing in {\em Rm1} of
(\ref{fig-st-isom}), and the {\em Me2}-move yield $S(G)=G^{-1}\,$,
which we thus have to impose in order for \stA\ to be non-trivial.

It is sometimes useful to have a ``co-tensorproduct'' on \stA, as in [KR] that
acts on a morphism by pushing a strand with labeling $a$ off itself 
in the plane  so that we get two parallel strands with labels $a'$ and $a''$.
Further, a blob with an element $x\in\A$ is duplicated, and the 
elements assigned to each of the new blobss are $x'_j$ and $x''_j\,$, where
$\Delta(x)=\sum_jx'_j\otimes x''_j\,$, see [KR]. In order for this to be 
compatible with {\em Rm1}, we have to impose  
the  somewhat stronger condition that $G$ is group-like, i.e.,
$\Delta(G)\,=\,G\otimes G\,$. A group-like element with (\ref{eq-G-S}) is
now exactly what is meant by a {\em balancing} of a general Hopf algebra.

In [KR] the  connection between \tgl\ and \stA\ is being made. We shall
review it here as a functor
\beq\lll{eq-tgl-stA}
{\sf Dec}_{\cal A}\;:\;\tgl\,\longrightarrow\,\stA\quad.
\eeq

It is clear that for this to exist we must be able to represent braid
group relations in \stA, too. Thus we have to assume that \A\ is in addition
quasi-triangular. This means that \A\ shall have an element $\cal R$ 
as in (\ref{eq-spec-elem}), which satisfies the Yang-Baxter equation.
Hence, if we define ${\sf Dec}_{\cal A}$ on a crossing in \tgl, as in
the right of diagram ({\ref{fig-Dec-cross}), tangles in \tgl, that
are equivalent under the third Reidemeister move, are mapped to equivalent
tangles in \stA.

\beq\lll{fig-Dec-cross}
\begin{picture}(470,50)

\put(20,45){\line(1,-1){40}}
\put(20,5){\line(1,1){17}}
\put(60,45){\line(-1,-1){17}}

\put(85,30){${\sf Dec}_{\cal A}$}
\put(80,25){\vector(1,0){37}}

\put(145,50){\line(1,-1){50}}
\put(145,0){\line(1,1){50}}
\put(155,40){\circle*{10}}
\put(160,44){$e_j$}
\put(185,40){\circle*{10}}
\put(194,38){$f_j$}

\put(250,45){\line(1,-1){17}}
\put(290,5){\line(-1,1){17}}
\put(250,5){\line(1,1){40}}

\put(308,30){${\sf Dec}_{\cal A}$}
\put(303,25){\vector(1,0){37}}

\put(375,50){\line(1,-1){50}}
\put(375,0){\line(1,1){50}}
\put(385,10){\circle*{10}}
\put(355,15){$S(e_j)$}
\put(415,10){\circle*{10}}
\put(422,13){$f_j$}
\end{picture}
\eeq

Now, we also can consider co-products of strands in \tgl . The condition that 
${\sf Dec}_{\cal A}\,$ maps this operation to the co-tensor product in
\stA\ translates with the assignment in (\ref{fig-Dec-cross}) into the
familiar triangle equations for $\cal R$ of a quasi-triangular 
Hopf algebra. 

Maxima and minima in \tgl\ shall be mapped by ${\sf Dec}_{\cal A}\,$
to their counterparts in \stA . Then the 
 second assignment in (\ref{fig-Dec-cross}) results from
checking the crossing-symmetry. Since 
$S\otimes 1({\cal R})\,=\,{\cal R}^{-1}\,$ we find that ${\sf Dec}_{\cal A}\,$
is also compatible with the second Reidemeister move in \tgl .

It remains to verify the framed version of the first Reidemeister move in
\tgl . To this end let us recall here that a quasi-triangular Hopf algebra has
a number of canonical elements, which can be constructed from the $R$-matrix, 
see [Dr]. The first of these is  $\,u\,=\,\sum_jS(f_j)e_j\,$. It does in fact
implement the square of the antipode but it is not group like. Another
element with these properties is $\hat u\,:=\,S(u)^{-1}\,$. It turns out
that $g:=u\hat u\,$ is group like, with $Ad(g)\,=\,S^4\,$. In other
words, although also quasi-triangular Hopf algebras do not have a canonical
balancing, at least the square of the balancing is canonical.
Indeed, it turns out that our functor is compatible with
the first Reidemeister move iff we have
\beq\lll{eq-bal-g}
g\quad=\quad G^2\qquad.
\eeq
A balancing with this property is what we mean by
 a {\em balancing of a quasi-triangular
Hopf algebra}. Sometimes this is equivalently described by the
central element $v\,=\,u G^{-1}\,$, which has to satisfy the 
equation 
\beq\lll{eq-mon-v}
{\cal M}\,:=\,{\cal R}'{\cal R}\,=\,v\otimes v\,\Delta(v^{-1})\;.
\eeq
It is easily worked out that it is assigned to a local $2\pi$-twist
of a  ribbon in \tgl . Since $v$ is central this extension of ${\sf Dec}_{\cal A}\,$
is compatible with moving a $2\pi$-twist through a crossing.
Imposing also $S(v)=v$  
allows us to move  such twists through extrema.

Let us summarize these findings in a more formal statement:

\btm[{[H],[KR]}]\lll{thm-Dec-KR}

For any balanced, quasi-triangular Hopf algebra, the functor 
${\sf Dec}_{\cal A}\,$
is well defined. 

(I.e., it factors into the equivalence classes.)
\etm

The functor ${\sf Dec}_{\cal A}\,$ can also 
be extended quite naturally to tangles with coupons,
given a  two-sided cointegral, $\lambda\in\A$. For a coupon, that has 
$m$ strands passing through we associate $m$ vertical strands, each eith one blob
 as in diagram
(\ref{fig-Dec-coup}). Here the elements assigned to the blobs are from the $m$-fold
coproduct $\Delta^{(m-1)}(\lambda)\,=
\,\sum_j\lambda^{(1)}_j\otimes\ldots\otimes\lambda^{(m)}_j\,$. This is again
compatible with the co-tensor product.

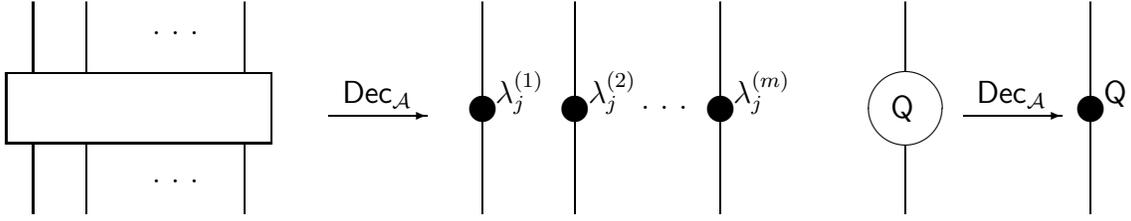
\beq\lll{fig-Dec-coup}
\begin{picture}(485,90)
\put(20,32){\framebox(100,26)}
\put(30,58){\line(0,1){27}}
\put(30,32){\line(0,-1){27}}
\put(50,58){\line(0,1){27}}
\put(50,32){\line(0,-1){27}}
\put(110,58){\line(0,1){27}}
\put(110,32){\line(0,-1){27}}
\put(75,73){$.\ .\ .$}
\put(75,17){$.\ .\ .$}

\put(147,47){${\sf Dec}_{\cal A}$}
\put(142,42){\vector(1,0){37}}

\put(200,5){\line(0,1){80}}
\put(200,45){\circle*{10}}
\put(205,47){$\lambda^{(1)}_j$}

\put(235,5){\line(0,1){80}}
\put(235,45){\circle*{10}}
\put(240,47){$\lambda^{(2)}_j$}

\put(260,45){$.\ .\ .$}

\put(290,5){\line(0,1){80}}
\put(290,45){\circle*{10}}
\put(295,47){$\lambda^{(m)}_j$}

\put(360,45){\circle{28}}
\put(360,59){\line(0,1){26}}
\put(360,31){\line(0,-1){26}}
\put(354.5,41.5){${\sf Q}$}

\put(387,47){${\sf Dec}_{\cal A}$}
\put(382,42){\vector(1,0){37}}

\put(430,5){\line(0,1){80}}
\put(430,45){\circle*{10}}
\put(435,47){{\sf Q}}
\end{picture}
\eeq

In general, we can extend ${\sf Dec}_{\cal A}$ to tangles with decorations 
by an assignment as in the right of diagram (\ref{fig-Dec-coup}) if the chips 
are labeled by central, $S$-invariant elements ${\sf Q}\in\A$.
The arguments are exactly the same as for $2\pi$-twists, which can be 
considered special examples of chips.
Thus  blobs with this labeling are of the same kind as the chips in
the previous paragraph. 
\medskip

Now, the coupons were in fact substitutes for 1-handle attachments. Let us therefore
impose in \tgl\ also as additional equivalences 
the $2\pi$-twist at a coupon and the coupon-slide from (\ref{fig-BL-mvs}). 
In order for ${\sf Dec}_{\cal A}$ to respect the  $2\pi$-twist, we have to 
require again  $v\lambda=\lambda\,$. The coupon-slide imposes, e.g.,  
${\cal M}(\lambda\otimes 1)=(\lambda\otimes 1)\,$, where $\cal M$ is as in
(\ref{eq-mon-v}). Given the $\varepsilon\otimes id({\cal R})\,=\,1$ and
$\varepsilon(v)=1$ this follows immediately from (\ref{eq-coint-class}).
Thus we can add the following:

\begin{rem}\lll{rem} 
${\sf Dec}_{\cal A}$ is also well defined for tangles with the coupons
of 1-handle attachments and other decorations that are labeled by
$S$-invariant  elements in ${\tt Z}_{\cal A}\,$.
\end{rem}

The algorithm to compute $\tau_{HKR}\,$ for a link, $\li$, now proceeds by 
unknotting its image in \stA\ into a bunch of disjoint circles
using the Reidemeister moves, with the singular version  {\em Rm2} from (\ref{fig-st-isom}).
By pushing the elements 
along each component together and multiplying them, we  arrive
at only one blob on each component of the link. $\tau_{HKR}$ is 
finally obtained by evaluation of the corresponding elements against $\mu^R\,$.

In a more precise language, that is also closer to the strategy by  which
$\tau_{RT}$ and $\tau_{L}\,$ are constructed, the result of the first part
of the algorithm, given by the rules in \stA, can be stated as follows:

\blm\lll{lm-Henn-alg}
Suppose $\li$ is a closed, projected  link with $N$ components, and 
$\li^{split}\,$ is a splitting as described in (\ref{fig-rib-split}).
Moreover, let \A\  be a balanced, quasi-triangular Hopf algebra.

Then there exists a unique element
\beq\lll{eq-split-elem}
\Ii\bigl(\li^{split},\,\A\bigr)\;=\;\sum_jA^{(1)}_j\!\otimes\ldots\otimes\! A^{(N)}_j\qquad\in\quad\A^{\otimes N}\;,
\eeq
such that we have the following equivalence in \stA\ :

\beq\lll{fig-li-Hred}
\begin{picture}(420,70)

\put(10,40){${\sf Dec}_{\cal A}\Bigl(\li^{split}\Bigr)$}
\put(100,40){$\simeq$}

\put(135,60){\rule{3.8in}{1mm}}

\put(150,60){\line(0,-1){15}}
\put(180,60){\line(0,-1){15}}
\put(165,45){\oval(30,30)[b]}
\put(180,45){\circle*{10}}
\put(184,39){$A^{(1)}_j$}

\put(220,60){\line(0,-1){15}}
\put(250,60){\line(0,-1){15}}
\put(235,45){\oval(30,30)[b]}
\put(250,45){\circle*{10}}
\put(254,39){$A^{(2)}_j$}

\put(290,67){${\bf R}_t$}
\put(300,43){$.\ \ .\ \ .$}

\put(360,60){\line(0,-1){15}}
\put(390,60){\line(0,-1){15}}
\put(375,45){\oval(30,30)[b]}
\put(390,45){\circle*{10}}
\put(394,39){$A^{(N)}_j$}

\end{picture}
\eeq
with corresponding labels at the top-line ${\bf R}_t\,$.
\elm

The existence is  quite obvious from our description of the algorithm.
A proof of uniqueness, i.e., the fact that the elements in 
(\ref{fig-li-Hred}) actually form a basis, 
 can be drawn from the antipode rule for pushing elements
through tangle pieces, together with the combinatorics given by isotopies and
the new {\em Rm1} rule. We shall not go into details here since 
Lemma~\ref{lm-Henn-alg} follows at once from the relations among categories
in the next paragraph.
\medskip

In this formalism the invariant is now obtained as follows:

\btm[{[H],[KR]}]\lll{thm-Henn-inv}
  Suppose $\li$ is an $N$-component link with splitting $\li^{split}\,$.
Let \A\ be a balanced, quasi-triangular Hopf algebra with 
$\Ii\bigl(\li^{split},\A\bigr)\,$ as in Lemma~\ref{lm-Henn-alg}.
For a tuple of states, $\mu_j\in\A^*\,$, with $j=1,\ldots,N\,$ 
define the number
\beq\lll{eq-II-Henn}
\tau_{HKR}\bigl(\li^{split},\,\A,\,\vec \mu\bigr)\;:=
\;\Bigl\lz \mu_1\otimes\ldots\otimes\mu_N\,,\,\Ii\bigl(\li^{split},\A\bigr)\Bigr\rz\;
=\;\sum_j\,\mu_1(A^{(1)}_j)\cdot\ldots\cdot\mu_N(A^{(N)}_j)\;.
\eeq
\ben
\item Suppose further that the states fulfill the following properties:
\ben
\item For all $j$, $\;\mu_j(xy)\,=\,\mu_j(S^2(y)x)\,$, or, equivalently,\\
 \ all $\mu_j$ are invariant under the 
$ad^*$-action as in Lemma~\ref{lm-coend-tann}.
\item For all $j$, $\;S^*(\mu_j)=\mu_j$, i.e., $\;\mu_j(S(x))\,=\,\mu_j(x)\,$. 
\een
Then $\tau_{HKR}$ is independent of the splitting.
\item Suppose $\mu_j=\mu^R\,$, i.e., the right integral of \A. Then
$\tau_{HKR}$ is invariant under 2-handle slides, and hence can be normalized
to a three-manifold invariant.
\een
\etm 

The claim made in (\ref{eq-S1S2=0}) of the introduction is now easily verified
by direct evaluation of the invariant on $S^1\times S^2\,$. A surgery 
diagram is given by an unframed unknot, $\li=\bigcirc^0$, so that we obtain 
$\Ii\bigl(\li^{split}\bigr)=1\,$, and hence
$$
\tau_{HKR}\bigl(S^1\times S^2\bigr)\;=\;\mu^R(1)\qquad.
$$
It follows from results in [LS] that $\mu^R(1)\neq 0$ is equivalent to
\A\ being cosemisimple. For modular Hopf algebras we have that the invariant
is also equal to $\varepsilon(\lambda)$, i.e., what we obtain if we present
$S^1\times S^2$  by a single coupon without strands instead of $\bigcirc^0\,$.
In particular, modular Hopf algebras are semisimple if and only if they are
cosemisimple.
\medskip

We shall see in the next paragraph that the constructions of
 $\tau_L$ and $\tau_{HKR}$ are equivalent so that Theorem~\ref{thm-Henn-inv}
can be inferred from Lemma~\ref{lm-int-inv}.
Nevertheless, let us sketch a few elements of the combinatorial  proof in
the category \stA :
\medskip

The ``$q$-cyclicity'' (or $ad^*$-invariance) in {\em (a)} ensures that for
two consecutive elements on a component $C$, we get the same result, whether
we attach the splitting ribbon $R_C$ before the pair or in between them.
A change of sides of the  attachment can be thought of an application of
the braided antipode $\Gamma\,$, which, by using $q$-cyclicity, can be reduced
to invariance under the ordinary antipode $S$. Also, changing an over-crossing
of a splitting ribbon into an under-crossing can be found from 
 the equation $S\otimes 1({\cal R})\,=\,{\cal R}^{-1}\,$, 
making  also use of {\em Me2}.
\medskip

The argument that translates the 2-handle slide into the defining equation 
of $\mu^R\,$ can be found in Kauffman's work, or it may be derived from
the analogous argument for categories used by Lyubashenko in [L],
combined with the equivalences shown here.
Alternatively, we may also try to check the 1-2-cancellation and modification
from (\ref{fig-BL-mvs}). They translate into the conditions $\mu(\lambda)=1$
and $(\id\!\otimes\!\mu)(\,\omega\,)\,=\,\lambda\,$ with $\omega\in\A^{\otimes 2}\,$
as in (\ref{eq-omega-A}). The latter was identified in 
Theorem~\ref{thm-int-modular} with the modularity property of \A.

\prg{3.4}{ Computing $\tau_L$ and $\tau_{RT}$ from $\tau_{HKR}\,$:
 Two Fiber Functors and a Central Element}
The strategy to prove the second part of Theorem~\ref{thm-main-rel}
 is to realize that the elements $\Ii\bigl(\li\bigr)\,$,
defined in Lemma~\ref{lm-int-inv} for the construction of $\tau_L$ and in
Lemma~\ref{lm-Henn-alg} for the computation of $\tau_{HKR}\,$, are identified
with the  same vectors by a suitable pair of fiber-functors. 
For large enough vector spaces this identification becomes one-to-one and
we can also relate their evaluations against integrals:
\medskip

For \C=\Am\  a fiber-functor is  of course given by  the forgetful functor ${\sf V}\,$ from 
(\ref{eq-tann-fct}). For \stA\ with the same algebra \A, we first have to
specify a coloring, {\tt Col}. This consists again of a choice of directions
for the strands of a tangle, and an assignment, $a\,\mapsto\,V_a\,$. Here
we think of $V_a$ as a vector space, on which we can also declare an \A-action.
In a more formal language, there  shall be an object $X_a\in\C$, such that
\beq\lll{eq-rel-col}
V_a\;=\;{\sf V}\bigl(X_a\bigr)\qquad\quad.
\eeq
The morals behind these formalities are that the objects in \C\ shall be
thought of as abstract objects as opposed to the vector spaces, $V_a$, 
and, furthermore,  we wish to permit operations on
the $V_a$, like transpositions of tensor factors, that make no sense for the objects in \C.

We then define a $\otimes$-functor 
\beq\lll{eq-fibfct-st}
{\sf W}^{\tt Col}\,:\;\stA\,\longrightarrow\,\vc
\eeq
on the generating morphisms of $\stA$. An object in \stA, i.e., a
 string of labels  $\{a_1,\ldots,a_k\}\,$, is mapped to the tensor product
$V_{a_1}^{\#}\otimes\ldots\otimes V_{a_k}^{\#}\,$ in analogy to the definition
of $I^{\tt Col}\,$ in (\ref{eq-fc-li-ab}), where $V^{\vee}$ is now simply the
dual space $V^*\,$.

Since ${\sf W}^{\tt Col}$ is supposed to respect $\otimes$,  we
map the  juxtapositions of two elementary diagrams to the respective 
tensor products of linear maps.  A functor is thus uniquely specified by 
the image of the elementary diagrams in (\ref{fig-sing-gen}).

The assignment for a single strand, carrying an element $x\in\A\,$, is 
simply the application of $x$ to $V$ if the direction is downward,
and of $S(x)^*$ to $V^*$ if it is upward.

To a crossing, $\,{\sf X}\,:\,\{a,b\}\to\{b,a\}\,$, we simply assign the 
transposition $T_{a,b}\,:\,V_a^{\#}\otimes V_b^{\#}\,\to\,V_b^{\#}\otimes V_a^{\#}\,$ of vector spaces.

The assignments to maxima and minima are as in (\ref{fig-rigid}). Here 
$ev_V$ and $coev_V$ are the canonical pairing and its inverse, and the flipped
morphisms are defined  by formula (\ref{eq-G-coev}).

It remains to show that ${\sf W}$ assigns the  same linear maps
to  equivalent composites:

\blm\lll{lm-Wcol-def}
The functor ${\sf W}^{\tt Col}\,$ from (\ref{eq-fibfct-st}) is well defined 
for every coloring {\tt Col}.
\elm

The relations that contain no blobs with elements from \A\ are immediate
from the fact that \vc\ is a symmetric, rigid tensor category with respect
to both choices of (co-)evaluations. Also {\em Me1} and {\em Me3} are obvious 
from the definition of the action of  an element $x\in\A$. 
The moves {\em Rm1} and {\em Me2} can be checked directly for both
choices of directions through the minimum, from (\ref{eq-G-coev})
and the fact that the condition $(a\otimes 1)\theta\,=\,(1\otimes S(a))\theta\,$
for all $a\in\A$ is equivalent to $\theta\in Inv_{\cal A}(V_1\otimes V_2)\,$.

We thus have fiber-functors on both \C\ and \stA. Recall, that we have also 
defined functors in  (\ref{eq-fc-li-ab}) and (\ref{eq-tgl-stA}) that
start at the category \tgl\ and end in the latter two categories. The 
intuitive correspondence, that guided Hennings and Reshetikhin in formulating
their rules, can thus be expressed  in  a more precise  language as a
commutative diagram of functors that describe different ways of 
assigning linear maps to purely topological tangles:

\btm\lll{thm-comm-fct}
Suppose we have chosen a colorings for $I^{\tt Col}\,$ from  
(\ref{eq-fc-li-ab}) and the one for ${\sf W}^{\tt Col}$ obtained from 
(\ref{eq-rel-col}). Then the following diagram of functors commutes.
\beq\lll{eq-comm-fct}
\bar{ccc}
\tgl&\TO{$I^{\tt Col}$}{60}&\Am\\
{\sf Dec}_{\cal A}\,\Biggl\downarrow\;\qquad&&
\quad\Biggl\downarrow\,{\sf V}\\
\stA&\TO{${\sf W}^{\tt Col}$}{60}&\vc
\ear
\eeq
\etm

  {\em Proof :} The proof  is 
easy since we only have to check the images of the 
generators.
It is clear that a label $a$ is mapped in both cases to $V_a$ or $V_a^*$
depending on direction. The action of an element that lives on a strand
in \stA\ on the vector-spaces coincides in both cases with the module
action of \A\ on $X_a$ and $X_a^{\vee}\,$.

Hence, the linear map 
$T_{a,b}{\cal R}_{V_a^{\#},V_b^{\#}}\,$, which we obtain by applying
${\sf W}^{\tt Col}\circ {\sf Dec}_{\cal A}$ to  
a braid on labels $a,b$, coincides with the definition of the commutativity
${\sf V}\bigl(\epsilon(X^{\#}_a,X^{\#}_b)\bigr)\,$ on specific \A-modules.
Consistency for coupons, $2\pi$-twists, and other central elements that may
occur as decoration in \tgl, follows by exactly the same arguments.

Finally, we remark that the (co-)evaluations produced by ${\sf W}^{\tt Col}\,$
are in fact invariant under the corresponding \A-actions. Thus they can be 
used to {\em define } a rigidity structure on \Am (and 
all such structures are  isomorphic).\hfill$\Box$
\medskip

This result allows us to make the first identification between the 
computations of $\tau_{HKR}$ and $\tau_L\,$.
Consider now a split link, $\li^{split}\,$, whose closure represents a manifold.
The first step to find $\tau_L$ was to compute the lift of the
transformations in  (\ref{eq-I-X-1}) to the coend. In the particular
situation of \C=\Am\ we have ${\sf V}(F)=\A^*$  so that the lift 
$\Ii_L(\li^{split})\,:\,F^{\otimes N}\to 1\,$ is realized
by {\sf V} in \vc\ by a linear map
$\,\Ii_L(\li^{split},\Am)'\,:\,(\A^*)^{\otimes N}\to 1$, or, equivalently,
by an element
\beq\lll{eq-Llift-tann}
\Ii_L(\li^{split},\,\Am)'\quad\in\quad\A^{\otimes N}\qquad.
\eeq

Recall from Lemma~\ref{lm-Henn-alg} that in the first step of computing
$\tau_{HKR}$ we also constructed an element  
$\Ii_H(\li^{split},\,\A)\,\in\,\A^{\otimes N}\,$ as in (\ref{eq-split-elem}).
Using  Theorem~\ref{thm-comm-fct} we can now  conclude equality:

\begin{cor}\lll{cor-I=I}
Let $\li^{split}$ be a split link, and \A\ a balanced, modular, quasi-triangular
Hopf algebra. Then
\beq\lll{eq-IL=IH}
\Ii_L(\li^{split},\,\Am)'\quad=\quad\Ii_H(\li^{split},\,\A)\qquad.
\eeq
\end{cor}

{\em Proof : } To begin with, let us choose a coloring, {\tt Col},
 of $\li^{split}\,$
given by a string of objects $X_1,\ldots,X_N\,$, and the corresponding
vector spaces $V_j\,:=\,{\sf V}(X_j)\,$. We now
determine the image of $\li^{split}$ in \vc\ in the both ways suggested
by Theorem~\ref{thm-comm-fct}, using the special elements in $\A^{\otimes N}\,$

Going through \C, $\li^{split}\,$ is mapped to the transformations $I^{(X_1,\ldots,X_N)}$
from (\ref{eq-I-X-1}), which can be written as a product of the lift
$\Ii\,$ to the coend and the transformations, $i_{X_j}\,$. Thus the image 
in \vc\ is also the product of $\Ii_L'\,=\,{\sf V}(\Ii)\,$ and the 
maps ${\sf V}(i_{X_j})\,:\,V^*_j\otimes V_j\to \A^*\,$, which are by 
Lemma~\ref{lm-coend-tann} given as the evaluation of matrix elements.
In summary, we have
\beq\lll{eq-abc}
\bar{rcl}
{\sf V}\Bigl(I^{(X_1,\ldots,X_N)}(\li^{split})\Bigr)\,&:
&\;V^*_1\otimes V_1\otimes\ldots\otimes V_N^*\otimes V_N\,\;\longrightarrow\quad\,{\bf C}\\
&&\\
\qquad l_1\otimes v_1\otimes\ldots \otimes l_N\otimes v_N\,&\mapsto&\,
\Bigl\lz l_1\!\otimes\ldots
\otimes\! l_N\,,\,\,\Ii_L(\li^{split})'(v_1\!\otimes\ldots\otimes\!v_N)\Bigr\rz
\ear
\eeq
for any choice of $l_j\in V^*_j\,$ and $v_j\in V_j\,$.

In the composite over \stA\ the link $\li^{split}$ is first mapped to a
series of $N$ arcs as in (\ref{fig-li-Hred}), each with a blob on it,
labeled by an element $A^{(k)}_j\,$ of \A.
It follows easily from the rules that, with directions going from
right to left, every one of these arcs is functored by ${\sf W}^{\tt Col}$ into
the linear map  $ev_V(1\otimes A)\,:\,V^*\otimes V\,\to\,{\bf C}\;
:l\otimes v\mapsto \lz l,\,Av\rz\,$. We thus obtain the same formula for
${\sf W}^{\tt Col}\Bigl({\sf Dec}_{\cal A}(\li^{split})\Bigr)\,$ as in
(\ref{eq-abc}), only with $\Ii_L$ replaced by $\Ii_H$ from 
(\ref{eq-split-elem}). Theorem~\ref{thm-comm-fct} therefore implies that
the equality (\ref{eq-IL=IH}) holds for all matrix elements of any choice of representations.
Since, e.g., the projective representation $Q$ from Paragraph 2.1 is faithful,
i.e., $\A\,\to\,End_{\bf C}(Q)\,$ is injective,
this implies the equality of the elements in $\A^{\otimes N}\,$ themselves.
\hfill$\Box$
\medskip

As we remarked in the end of Paragraph 2.5, the categorical integral, $\mu$,
of $F$ is identical with the classical, right integral, $\mu^R$, of \A.
We thus have the equality,
$$
\tau_L\;=\;\bigl\lz \mu^{\otimes N} , \Ii_L\bigr\rz\;=
\;\bigl\lz (\mu^R)^{\otimes N}, \Ii_H\bigr\rz\,=\,\tau_H\quad,
$$
which completes the proof of Theorem~\ref{thm-main-rel}.\hfill $\Box$
\medskip

Recall, that the element  from (\ref{eq-Llift-tann}) actually defined an
intertwiner, $\Ii_L\,:\,(\A^*)^{\otimes N}\to 1$, with the coadjoint
action defined on $\A^*$. We thus have as another, immediate
 application of the correspondence in Corollary~\ref{cor-I=I} the following:

\begin{cor}\lll{cor-IHinvar}
The element $\Ii_H\bigl(\li^{split},\,\A\bigr)\,\in\,{\cal A}^{\otimes N}\,$,  
computed from the rules in \stA, is invariant under the \A-action, defined
by the $N$-fold-tensor product of the adjoint representations.

In particular, for a split knot ${\cal K}^{split}$  ($N=1$) we obtain a 
central elements, 
$\Ii_H\bigl({\cal K}^{split},\,\A\bigr)\,\in\,{\tt Z}_{\cal A}\,$,
from Hennings' rules.
\end{cor}

Since a split knot is the same as a one-component tangle the last 
assertion is exactly about the centrality of the elements computed
in [Re] and the problem that L. Kauffman addressed in his lecture.
\bigskip 

In Theorem~\ref{thm-RT-LQ} we learned about a modification of the
construction of $\tau_L$ that yields $\tau_{RT}$, even if \C\ is not
semisimple. Let us use the equivalence between $\tau_{L}$ and $\tau_{HKR}$
in order to derive an analogous rule in the Hennings picture.

In the  modification for general categories, we considered
links, $\li^{\sf Q}\,$, where
we have decorated each component with a chip labeled  by a natural
transformation of the identity, {\sf Q}.
By Lemma~\ref{lm-coend-1} this corresponds for $\C=\Am\,$  to the application of 
a central element
$$ 
{\sf Q}\in {\tt Z}_{\cal A}\qquad.
$$
Moreover, condition (\ref{eq-jjj}) translates into $S$-invariance of {\sf Q}.
By Remark~\ref{rem} and the rule indicated in (\ref{fig-Q-ins}) 
we can extend the functors in Theorem~\ref{thm-comm-fct} to tangles with 
decorations. If the corresponding natural transformation for $I^{\tt Col}\,$
is defined by  application of the central element used for insertions 
in ${\sf Dec}_{\cal A}\,$ it is clear that (\ref{eq-comm-fct}) still
commutes.
\medskip

This implies that the invariant $\tau_L(\li^{\sf Q})\,$  of a decorated link
 can be computed equivalently from $\tau_H(\li^{\sf Q})\,$, generalizing 
Theorem~\ref{thm-main-rel}. Before we use this to reformulate the result in Theorem~\ref{thm-RT-LQ},
 let us give a definition of 
the central element {\sf Q} in the setting of an ordinary Hopf algebra:

The definition from (\ref{def-Q-bms}) translates to 
$\mu^{semis}\,=\,\mu^R\rcoa {\sf Q}\,$, where we used  the
sum of quantum-traces from (\ref{eq-mu-Css}). By virtue of (\ref{eq-trq-tann}),
we can characterize {\sf Q} uniquely by the following concrete
formula
\beq\lll{eq-defQ-H}
\mu^R\bigl({\sf Q}\cdot y\bigr)\;=\;{\cal D}^{-1}
\sum_{j\in{\cal J}_0}\,d(j)\,tr^{can}_{V_j}\bigl(G\cdot y\bigr)\qquad\forall\;y\in\A.
\eeq
The sum runs over a complete set of inequivalent, irreducible representations, 
$V_j$,  of \A, $G$ is the balancing, $tr^{can}$ is the canonical trace, and 
$d(j)=tr^{can}_{V_j}(G)\,$. 
Since the $tr_V^{can}\rcoa G\,$ and $\mu^R$ are all $q$-cyclic, and since $\mu^R$ is 
non-degenerate, it follows that any solution, {\sf Q}, of (\ref{eq-defQ-H}) is central.
Because of $d(j)=d(j^{\vee})\,$ and the properties of the quantum-traces,
 we have that $\mu^{semis}$ fulfills
$S^*(\mu^{semis})\,=\,G^{-2}\lcoa\mu^{semis}\,$. It follows from
 results  in  [Rd], see also 
[Ke2], that the same equation holds for $\mu^R\,$.
This can be used to show that {\sf Q} is  invariant under the antipode $S$.

\begin{cor}\lll{cor-HQ=RT}
Suppose \A\ is a balanced, quasi-triangular Hopf algebra  and 
${\sf Q}\in{\tt Z}_{\cal A}$ is as in (\ref{eq-defQ-H}). Then
$$
\tau_{RT}(\li,\A)\;=\;\tau_{HKR}(\li^{\sf Q},\A)\quad,
$$
where the additional rule for evaluating the decorated tangle,
is to insert the element {\sf Q} in each component of the link in
the \stA-picture.
\end{cor}

As opposed to expressing the relation between $\tau_{RT}$ and $\tau_{HKR}\,$
in a form analogous to (\ref{eq-rrr}), 
the picture given above lends itself much better to explain 
relations of the corresponding TQFT's, as we shall indicate in the next 
paragraph.
\medskip 

Quite obviously we have {\sf Q}=1 if \A\ is semisimple, since the formula 
$\mu^{semis}\,$ was computed for a general semisimple category. 
The converse is also true. If \A\ is not semisimple it follows that 
$\lambda^2=0$ so that
we must have $tr_{V}^{can}(G\cdot\lambda)=0\,$. Thus if ${\sf Q}=1$, i.e.,
if $\mu^R$ is given by a sum of traces, this implies $\mu^R(\lambda)=0\,$.
But that is not possible by the Lemma~\ref{lm-beta-inv}.

In fact a stronger statement holds, namely that {\sf Q} is always nilpotent:

\blm\lll{lm-Q2=0}Suppose \A\ is a non-semisimple, balanced, finite dimensional,
quasi-triangular Hopf algebra, and {\sf Q} is as above. Then
\beq\lll{eq-Q2=0}
{\sf Q}\cdot{\sf Q}\quad=\quad 0\qquad\;.
\eeq
\elm

{\em Proof: } Using the $S$-invariance of {\sf Q} and a relation from [Rd],
similar to the one from Lemma~\ref{lm-beta-inv}, we find that 
${\sf Q}\,=\,\mu^{semis}\!\otimes\!\id\bigl(\Delta(\lambda)\bigr)\,$, which can be
rewritten as a weighted sum over the elements 
${\sf Q}_j\,=\,\sum_{\nu}tr^q_j(\lambda_{\nu}')\lambda_{\nu}''\,$ for $j\in{\cal J}_0^*\,$.
Using the identity 
$\, (1\!\otimes\!S(a))\Delta(\lambda)\,=\,(a\!\otimes\!1)\Delta(\lambda)\,$,
and the multiplicativeness of the quantum-traces, we find for these elements
$$
S({\sf Q}_j){\sf Q}_k\,=\quad\sum_{\nu_1,\nu_2}\,tr^q_{j\otimes k}
\bigl(\lambda_{\nu_1}'\otimes\lambda_{\nu_2}'\bigr)\,
S(\lambda''_{\nu_1})\lambda''_{\nu_2}\;=\quad\sum_{\nu}\,
tr^q_{j\otimes k}\bigl(\Delta(\lambda)(1\!\otimes\!\lambda_{\nu}')\bigr)\,\lambda_{\nu}''\qquad.
$$
Now, the image of $\Delta(\lambda)$ on $V_j\otimes V_k\,$ is in the invariance,
which is zero, unless $j\cong k^{\vee}\,$. Moreover, if $d(k)\neq 0\,$ the 
invariance is a direct summand of $V_k^{\vee}\otimes V_k\,$ so that 
$\Delta(\lambda)\,$ acts as $\varepsilon(\lambda)=0\,$. It follows that
$S({\sf Q}){\sf Q}={\sf Q}^2=0\,$. \hfill $\Box$
\smallskip

We shall given a more explicit
 form of {\sf Q} in the case of $U_q(s\ell_2)\,$ in the next paragraph. 

\prg{3.5}{ Generalization of Hennings' Rules to  TQFT's}
The extend, to which it is possible to construct TQFT's from non-semisimple
categories, has been explained in Paragraph 1.5. In particular, 
Theorem~\ref{thm-conn-TQFT} together with the relation between $\tau_L$ and
$\tau_{HKR}\,$ shows that there exists a TQFT
for connected surfaces that specializes to $\tau_{HKR}\,$. 
The aim of this paragraph is to outline a set of rules, which generalize those
of \stA, and which  allow us to compute the linear maps associated to a 
cobordism.     

>From these rules we derive representations of the 
mapping class groups of the torus, that had been discovered before in [FLM]. 
In the work of Lyubashenko and Majid the formulas were obtained
using rather similar ideas as the ones presented here. 
Quite interestingly, they also appeared in
works of Felder et al. on  local systems that arise in the study of integral
representations of conformal blocks.

We shall review  the detailed structure of these representation
, that was obtained in [Ke2] for the case of $U_q(sl_2)\,$, 
and derive from this a relation between $\tau_{HKR}$ and $\tau_{RT}\,$ for
lens spaces. 
We also show how the special element {\sf Q} can be used to directly relate
the (inequivalent)  representations  arising in both pictures.
\medskip

Let us sketch first a few ingredients of the combinatorial representation
of the category $Cob_3^{conn}$ in terms of tangles, as proposed in [Ke5].
 Here  $Cob_3^{conn}$ is the category of connected three-manifolds, $M$,
 (or more precisely a central extension by $\Omega_4$ thereof) 
that cobord one connected surface to another connected surface, i.e.,
$\partial M\,\cong\,-\Sigma_{g_1}\amalg\Sigma_{g_2}$ so that 
$M\,:\,[g_1] \to [g_2]\,$. 
(compact and orientable always assumed).

A tangle ${\cal T}_M$ in $\tgl$ that represents $M$ will have $2g_1$ endpoints
at the top-line ${\bf R}_t$ and $2g_2$ endpoints at the bottom-line
${\bf R}_b\,$. In one version of the presentation-calculus  
 every component of an admissible tangle must be  one of the following
three types:
\ben
\item A {\em top-ribbon}, that starts at the $2j$-th position at ${\bf R}_t$
 and returns to the $2j-1$-st position at ${\bf R}_t\,$ ,
\item a {\em bottom-ribbon}, that starts at the $2j$-th position at ${\bf R}_b$
 and returns to the $2j-1$-st position at ${\bf R}_b\,$, or
\item a {\em closed} ribbon, that is disjoint from the boundaries of
the diagram.
\een

Hence there are exactly $g_1$ top-ribbons, $g_2$ bottom-ribbons, and any number, $N$,
of closed ribbons. Moreover, the label structure of ${\cal T}_M$ as a morphism in
\tgl\ is as follows:
$$
{\cal T}_M\;:\;\{a_1,a_1,\ldots,a_{g_1},a_{g_1}\}\;\longrightarrow
\;\{b_1,b_1,\ldots,b_{g_2},b_{g_2}\}\qquad.
$$
It will be convenient to use a slightly larger space, where we admit
in addition pairs of {\em through-ribbons} that can be either regular or crossed.
A pair of such ribbons is associated to a pair of labels, $(a_j, b_k)$, 
one from the top, one from the bottom. In a {\em regular} pair of through-strands
one strand connects the first $a_j$-label at ${\bf R}_t\,$
to the first $b_k$-label at ${\bf R}_b\,$ and the 
second starts and ends at the other two labels. For a {\em crossed} through-pair
we have one ribbon starting at the first $a_j$-label and ending at the second
$b_k$-label and the other ribbon connects again the remaining two labels.
\medskip

The moves for surgery presentations of closed manifolds shall also apply to regions of tangles in 
  \tgl\ that stay away from the boundaries of the diagram. In addition
we have two moves in the vicinity of the boundary. The first is what we call
the $\sigma$-move, which is  depicted below in the version with coupons:

\beq\lll{fig-sigma-mv}
\begin{picture}(400,130)

\put(65,115){$a_j$}
\put(105,115){$a_j$}
\put(70,105){\line(0,-1){85}}
\put(110,105){\line(0,-1){85}}

\put(52,105){\rule{1in}{1mm}}

\put(145,55){\vector(1,0){50}}
\put(150,60){$\sigma$-move}

\put(235,115){$a_j$}
\put(275,115){$a_j$}

\put(220,105){\rule{2.2in}{1mm}}
\put(240,105){\line(0,-1){55}}
\put(280,105){\line(0,-1){35}}
\put(260,50){\oval(40,40)[b]}

\put(270,50){\framebox(40,20)}

\put(320,70){\oval(40,40)[t]}
\put(340,70){\line(0,-1){50}}
\put(300,50){\line(0,-1){30}}

\end{picture}
\eeq

Invariance of the general construction in [KL] under the $\sigma$-move 
translates directly to  inversion-formulae as the one in 
Lemma~\ref{lm-beta-inv}. Under strong enough modularity conditions,
or the canonical choices for $\mu$ and $\lambda$,  
the invariance with respect to $\sigma$-moves  follows therefore
immediately.
\medskip

We also have the $\tau$-move, which is basically a special isotopy
of tangles over the sphere through the point at infinity. 
It allows us to push a strand, that is parallel and 
close to ${\bf R}_t$, and that crosses all strands emerging from the top-line,
to the back of everything, as it is indicated in the next diagram:

\beq\lll{fig-tau-mv}
\begin{picture}(400,100)
\put(33,93){$a_1$}
\put(57,93){$a_1$}
\put(35,85){\line(0,-1){58}}
\put(60,85){\line(0,-1){50}}

\put(35,23){\line(0,-1){15}}
\put(60,31){\line(0,-1){23}}

\put(78,70){$.\ \ .\ \ .$}

\put(133,93){$a_g$}
\put(135,85){\line(0,-1){25}}
\put(135,56){\line(0,-1){48}}

\put(15,85){\rule{2.2in}{1mm}}
\put(20,20){\line(3,1){150}}

\put(190,45){\vector(1,0){50}}
\put(195,50){$\tau$-move}

\put(283,93){$a_1$}
\put(307,93){$a_1$}
\put(285,85){\line(0,-1){77}}
\put(310,85){\line(0,-1){77}}

\put(330,30){$.\ \ .\ \ .$}

\put(383,93){$a_g$}
\put(385,85){\line(0,-1){77}}

\put(265,85){\rule{2.2in}{1mm}}
\put(270,20){\line(3,1){11}}
\put(288,26){\line(3,1){19}}
\put(387,59){\line(3,1){33}}

\end{picture}
\eeq

As a matter of fact, the omission of this move would give us a tangle
category that is equivalent to $\,Cob^{conn}_3(1)\,$, i.e., cobordisms between surfaces
with one hole.
\medskip

In [KL] we now construct a linear map for every  cobordism $M$, represented by a tangle
${\cal T}_M\;$. The first step is, again, to consider  a {\em split}
  tangle ${\cal T}^{split}\,$.
In this generalization we attach one   splitting-ribbon to every
closed and every bottom-ribbon. The top- and through-ribbons are not 
split.

Hence, if $\cal T$ has $N$ closed ribbons, $M$ pairs of through-ribbons,
$g_1-M$ top-ribbons, and $g_2-M$ bottom-ribbons, then the split tangle
${\cal T}^{split}\,$ will have $g_2$ pairs of through-ribbons, $(g_1-M)+N$
top-ribbons, no bottom-ribbons, and no closed ribbons. If we assume for 
simplicity that the labels of the through-ribbons in $\cal T$ are all to the 
left of the labels of the top-ribbons, we obtain a schematic
label structure:
$$
{\cal T}^{split}\,:\;
\Bigl\{L^M_{thr.},\,L^{g_1-M}_{top},\,L^{g_2-M}_{bot.},\,L^{N}_{cl.}\,\Bigr\}\;\longrightarrow\;\Bigl\{L^M_{thr.},\,L^{g_2-M}_{bot.}\,\Bigr\}\qquad,
$$
where $L^n_{\#}\,$ has $2n$-labels associated to components of $\cal T$ with
characteristic $\#$.
\medskip

In the next step we construct again  morphisms in \C\ 
by chosing  colorings, $a\to X_a$, and check their 
(di-)naturality properties. They allow us to lift everything to a single 
morphism
$$
\Ii\bigl({\cal T}^{split}\bigr)\;:\;F^{\otimes M}\otimes F^{\otimes(g_1- M)}\otimes
F^{\otimes (g_2-M)}\otimes F^{\otimes N}\quad\,\longrightarrow\,\qquad F^{\otimes g_2}\quad.
$$
>From this we obtain a morphism $\Ii_L\bigl({\cal T}\bigr)\,:\,F^{g_1}\to 
F^{g_2}\,$ by multiplying with $\id_{F^{g_1}}\otimes \mu^{\otimes (g_2-M+N)}\,$
from the right. If the tangle represented a cobordism in $\,Cob_3(1)^{conn}$, this would already
be the construction of the morphism in \C, that is associated by the TQFT-functor
${\cal V}_1\,$.
\medskip

In order to get a TQFT for closed surfaces, we also have to take care of 
the $\tau$-move. As explained in Paragraph 1.6 from general arguments, this can 
be done by composing ${\cal V}_1$ with the invariance. The vector space 
associated to a surface $\Sigma_g$ is then $Inv_{\cal C}(F^{\otimes g})\,$.

Our aim here is to find a functor that is explicitly described by 
Hopf algebra elements. We shall therefore construct instead the 
{\em contravariant} functor given by application of the coinvariance:
$$
{\cal V}^*_0\;=\;Coinv\circ{\cal V}_1
$$
The functor, ${\cal V}^*_0$ respects now the $\tau$-move  in the same way as ${\cal V}_0$ does.
 The vector space associated by ${\cal V}_0^*\,$ to $\Sigma_g\,$  is  
the $ad$-invariant 
 subspace $(\A^{\otimes g})_{ad}\,$ of $\A^{\otimes g}$. 
The image of a tangle under ${\cal V}^*_0$ is then  obtained by restricting the 
$\Ii\bigl({\cal T}\bigr)^*\,:\,\A^{\otimes g_2}\to\A^{\otimes g_1}$ to the 
$ad$-invaraiant subspaces.
The rules to compute the $\Ii\bigl({\cal T}\bigr)^*$ are the following:
\medskip

We start from a split, ${\cal T}^{split}$, of the representing tangle as above 
and map it to an element in \stA\ by applying ${\sf Dec}_{\cal A}$ from
(\ref{eq-tgl-stA}). It is not hard to see that by invoking the relations of 
\stA\ we can bring the singular tangle in the form described in the 
diagram below:

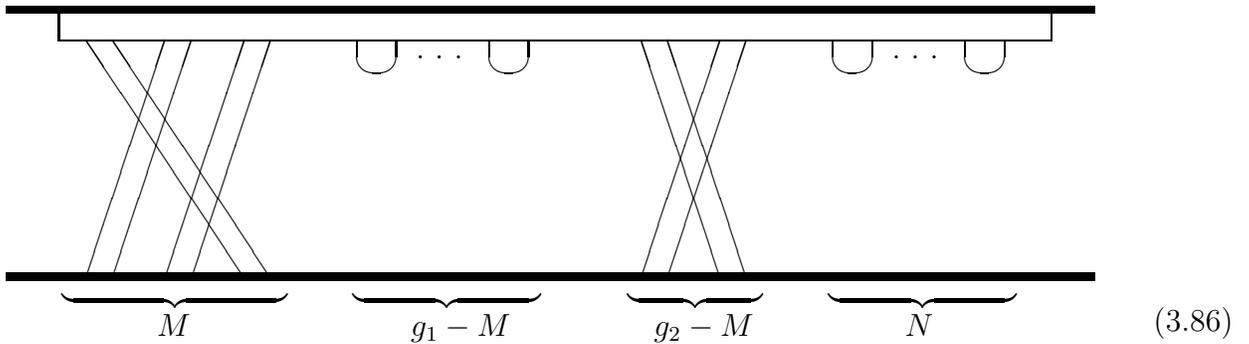
\beq\lll{fig-sd-norm}
\begin{picture}(420,130)

\put(10,19){\rule{5.7in}{1mm}}

\put(10,120){\rule{5.7in}{1mm}}

\put(30,110){\framebox(375,10)}

\put(40,20){\line(1,3) {30}}
\put(50,20){\line(1,3) {30} }

\put(70,20){\line(1,3) {30}}
\put(80,20){\line(1,3) {30} }

\put(100,20){\line(-2,3) {60}}
\put(110,20){\line(-2,3) {60}}

\put(31,15){$\underbrace{\hspace*{3cm}}$}
\put(67,-2){$M$}

\put(141,15){$\underbrace{\hspace*{2.5cm}}$}
\put(163,-2){$g_1-M$}

\put(245,15){$\underbrace{\hspace*{1.8cm}}$}
\put(255,-2){$g_2-M$}

\put(321,15){$\underbrace{\hspace*{2.5cm}}$}
\put(350,-2){$N$}

\put(150,110){\oval(15,25) [b]}
\put(165,103){$.\ .\ .$}
\put(200,110){\oval(15,25) [b]}

\put(250,20){\line(1,3) {30}}
\put(260,20){\line(1,3) {30}}

\put(280,20){\line(-1,3) {30}}
\put(290,20){\line(-1,3) {30}}

\put(330,110){\oval(15,25) [b]}
\put(345,103){$.\ .\ .$}
\put(380,110){\oval(15,25) [b]}

\end{picture}
\eeq

Here the strands of each of the $g_2$ pairs of through-strands 
are straight and parallel lines. Hence the diagram is uniquely characterized
by a permutation $\pi\,\in\,S_M\times S_{g_2-M}\,$ of label-pairs.

In the small strip at the top line we are allowed to have insertions of 
\A-elements and crossings of a pair of crossed through-strands.
They are depicted below:

\beq\lll{fig-ds-ins}
\begin{picture}(380,40)
\put(20,5){\line(0,1){30}}
\put(40,5){\line(0,1){30}}
\put(20,20){\circle*{10}}
\put(40,20){\circle*{10}}
\put(7,23){$p$}
\put(48,23){$q$}

\put(60,18){$\;\;\sim\quad\, x\mapsto S(p)xq\quad,$}

\put(200,5){\line(1,1){30}}
\put(230,5){\line(-1,1){30}}

\put(255,18){$\!\sim\qquad x\mapsto S(x)G\quad.$}

\end{picture}
\eeq

We consider possible insertions at the crossing as a composite of
the two diagrams in (\ref{fig-ds-ins}). 
For the minima we can use relations in \stA\ to 
assume a form as in diagram (\ref{fig-li-Hred}). 
Note that strands with different colorings, but also
the strands of a pair of regular through-strands, do not cross in the upper part
in this normal form. 
\medskip

A map $\A^{\otimes g_2}\to\A^{\otimes (M+(g_1-M)+(g_2-M)+N)}\,$  
may now be obtained by thinking of a  pair of  through-strands
as a  single strand, and a diagrams as in (\ref{fig-ds-ins}) as one
insertion. The translation into  maps on tensor products of \A\
 is then as for the functor ${\sf W}^{\tt Col}\,$, but now with morphisms going
from bottom to top. 

Specifically,
we first apply to an element in $\A^{\otimes g_2}$ the permutation $\pi$ 
indicated in (\ref{fig-sd-norm}). 
The morphism for a fixed summation index
of the insertions is further described by applying the maps associated in
(\ref{fig-ds-ins}) to each tensor factor (Note, that $G$ on the right side is the balancing so that the operation associated to a crossing is in fact an involution).
Finally, we also insert the elements
in $\A^{\otimes(g_1-M)}$ and $\A^{\otimes N }$, defined as in 
Lemma~\ref{lm-Henn-alg}, into the respective places in the tensor product,
and take the sum of the morphisms over all labels.
\medskip

This yields in fact a contravariant representation of $Cob_3(1)^{conn}\,$.
Generalizing the arguments of the previous chapter and using the results from
[KL] we thus find the following:

\btm\lll{thm-Henn-TQFT}
The algorithm described above yields a well defined, 
contravariant functor
$$
{\cal V}_0^*\,:\;Cob_3^{conn}\,\longrightarrow\,\vc
$$
with 
${\cal V}_0^*\bigl(\Sigma_g\bigr)\,\cong\,\Bigl(\A^{\otimes g}\Bigr)_{ad}\,$.

It is dual to the functor constructed in [KL], and specializes to $\tau_{HKR}$
for closed manifolds.
\etm

In the remainder of this chapter let us consider as an application the
case of the mapping class group $\pi_0\bigl(Diff(\Sigma_1)^+\bigr)\,$
of the torus. The standard generators, denoted by $T$ and $S$, are realized  as invertible 
cobordisms $\Sigma_1\to\Sigma_1\,$ by the following tangles:

\beq\lll{fig-ST-tgl}
\begin{picture}(260,60)

\put(-7,30){$T\,=$}

\put(30,5){\line(0,1){50}}
\put(60,5){\line(0,1){15}}
\put(60,40){\line(0,1){15}}
\put(50,20){\framebox(20,20){$2\pi$}}

\put(135,30){$S\,=$}

\put(197,40){\oval(40,40)[bl]}
\put(203,40){\oval(40,40)[br]}
\put(220,20){\oval(40,40)[tl]}
\put(226,20){\oval(40,40)[tr]}

\put(177,40){\line(0,1){15}}
\put(223,40){\line(0,1){15}}
\put(200,20){\line(0,-1){15}}
\put(246,20){\line(0,-1){15}}

\end{picture}\hfill .
\eeq

By a straightforward application of the above rules we find
 that the maps $T^*,S^*\,;\,\A\to\A\,$,
that are associated by our algorithm to these tangles, have the form:
\beq\lll{eq-ST-1}
T^*(x)\,=\,xv\,=\,vx\;,\qquad{\rm and}\qquad S^*(x)\,=\,\sum_j\,\omega^{(1)}_j\,\mu^R
\bigl(x\omega^{(2)}_j\bigr)\qquad,
\eeq
where $v$ is the central balancing element from (\ref{eq-mon-v}) and
 $\omega\,=\,\sum_j\omega^{(1)}_j\!\otimes\!\omega^{(2)}_j\,$ is as in
(\ref{eq-omega-A}). This representation of 
$\pi_0\bigl(Diff(\Sigma_{1,1})^+\bigr)\,$ is equivalent to the ones found in [FLM]

As we remarked following its definition,  the 
form $\omega$ has  particularly nice properties if $\A=D({\cal B})\,$ is a 
Drinfel'd double. In [Ke2] we find in this case a factorization, given
by the following commutative diagram. In the last chapter we will attempt to  
relate this  to a conjecture on relations between three-manifold invariants.

\beq\lll{eq-S-fact}
\bar{ccc}\B\otimes \Bo
&
\raise 1.35ex\hbox{$\underline{ \beta_l\circ S^{-1}\,\otimes \,S\circ\overline{\beta_l}}$}\mkern -6mu\to
&
\Bo\otimes\B
\\
\Bigg\downarrow {\bf \cdot}
&&
\Bigg\downarrow {\bf \cdot}
\\
D(\B)
&
\raise .83ex\hbox{$\underline{ \,\,\quad\qquad { S^*  } \,\,\quad\qquad}$}\mkern -6mu\to
&
D(\B)
\ear 
\eeq

Here, the vertical arrows are the isomorphisms given by multiplication in $D(\B)\,$.
The canonical maps in the upper horizontal arrow are defined by $\beta_l(x)(y)\,=\mu^L(yx)\,$,
and $\overline{\beta_l}(\rho)\,=\,\rho\rcoa\lambda^L\,$, where $\mu^L$ and $\lambda^L\,$ are the
left integral and cointegral of $\cal B\,$.
\medskip

In order to obtain maps  for the closed torus we have to restrict $T^*$ and $S^*$
to the center ${\tt Z}_{\cal A}\,$  of \A.
The projective representations of $SL(2,{\bf Z})\,$ we get in this way 
from $U_q(s\ell_N)\,$  at a $p$-the root of
unity are closely related to the representations we obtain from the Fourier transformations
on the group ${\bf Z}/p\,$ and its powers. If $p$ is an odd prime this representation has
two irreducible factors, given by the eigenspaces of the inversion-map:
\beq\lll{eq-Zp-dec}
C({\bf Z}/p)\;=\;{\cal Y}_{\frac {p+1} 2}^+\,\oplus\,{\cal Y}_{\frac {p-1} 2}^-\qquad,
\eeq
where the subscripts indicate the dimensions. 

The next theorem has been conjectured (with a clarification due to V. Jones) in [Ke2], and
the gap in the arguments can be  filled using the observation following it.

\btm\lll{thm-sl2-dec} Suppose $\A=U_q(s\ell_2)\,$ with $q$ a primitive $p$-th root of unity and $p$ an odd prime.
Then the (projective) representation of $SL(2,{\bf Z})\,$ on the center of \A, 
which is obtained  from the TQFT described above, has the following structure:
$$
{\tt Z}_{\cal A}\;=\;{\cal Y}_{\frac {p+1} 2}^+\;\,\oplus\;\,\Pi\otimes{\cal Y}_{\frac {p-1} 2}^-\qquad.
$$
In this decomposition the $\cal Y$'s are the finite representations from (\ref{eq-Zp-dec}),
 and $\Pi$ is the restriction of
the two-dimensional fundamental representation of $SL(2,{\bf R})\,$.
\etm

Although this result was obtained mainly by involved computations, it is possible 
to give some algebraic meanings to the constituents of ${\tt Z}_{\cal A}$.  
To begin with, we know that in the pairing of
invariance and coinvariance of an object in a non-semisimple category is 
usually degenerate. The null-space of this pairing is certainly an invariant subspace of
the representation. 

In our example this subspace turns out to be precisely the  summand ${\cal Y}_{\frac {p+1} 2}^+\,$.
A more concrete characterization results from the fact that  if we consider \A\ as an \A-module with
adjoint action the invariance is ${\tt Z}_{\cal A}$ and the coinvariance consists
of the $q$-cyclic states.  ${\cal Y}_{\frac {p+1} 2}^+\,$ is thus given by 
the central elements on which every $q$-cyclic state vanishes.

It is also known that the representations we
obtain from the semisimple picture are equivalent to the ${\cal Y}_{\frac {p-1} 2}$'s.
\medskip

By its very property the null-space does not contribute to the invariants of 
closed manifolds so that it suffices to consider the respective factor 
representation. In the situation of Theorem~\ref{thm-sl2-dec} this is also 
 a direct summand, and the generators are explicitly represented by 
the following matrices:

\beq\lll{eq-ST-prod}
T_H\;=\;\left[\bar{cc} 1&0\\1&1\ear\right]\!\otimes\!T_{RT}\qquad\qquad\;
S_H\;=\;\left[\bar{cc} 0&-1\\1&0\ear\right]\!\otimes\!S_{RT}
\eeq

Here, the subscript $_{RT}$ indicates, as usual, 
that the matrices are those that we obtain 
from the TQFT associated to the semisimple quotient of \Am. 
\medskip

Recall from Corallary~\ref{cor-HQ=RT} that  the computation of $\tau_{RT}\,$
can be done in the Hennings picture if we insert at each component of the link
the special element {\sf Q}, that was defined in (\ref{eq-defQ-H}).

Now, if a link is given by the composition  of many admissible tangles, with ribbon 
types defined as in the beginning of this paragraph, then there is a one-to-one 
correspondence between the closed and bottom ribbons of the composites and 
the components of the link. Thus, if we define a modified Hennings procedure
for TQFT's by inserting one {\sf Q} at every bottom or closed ribbon in a given tangle, 
then the product of  the assigned linear maps gives us $\tau_{RT}\,$.

It is, however, not true that this assignment defines a TQFT, since there are
degeneracies so that we usually do not preserve invariance under the 
$\sigma$-move. Often,  this can be  salvaged by taking proper quotients of
the vector spaces of the original TQFT:
\medskip

We know from Lemma~\ref{lm-Q2=0} that  {\sf Q} is nilpotent. In fact, 
multiplication by {\sf Q} it is represented for our special example
(up to a scalar) on the second 
summand in Theorem~\ref{thm-sl2-dec} by
$$
{\sf Q} \;=\;\left[\bar{cc} 0&0\\1&0\ear\right]\!\otimes\!\id\qquad\quad.
$$

If we wish to apply our insertion procedure to the $T$-generator,
we immediately see that nothing has to be done, since we have only through-strands,
i.e., $T^{\sf Q}_H\,=\,T_H\,$. However, the $S$ generator contains a bottom
ribbon: 
$$
S^{\sf Q}_H\;=\;S_H\cdot {\sf Q}\;=\;\left[\bar{cc} -1&0\\0&0\ear\right]\!
\otimes\!S_{RT}\qquad\quad.
$$
It is quite obvious that the matrices $T^{\sf Q}_H$ and $S^{\sf Q}_H$ do not
define  a representation of  $SL(2,{\bf Z})\,$ anymore. Also, while $S_H$ and
$T_H$ define an {\em irreducible} representations, the matrices 
$T^{\sf Q}_H$ and $S^{\sf Q}_H$ have $ker\bigl(S^{\sf Q}_H\bigr)=ker({\sf Q})\,$
as a common subspace. If we divide by this subspace we get again a representation 
of $SL(2,{\bf Z})\,$, namely exactly the one defined by $S_{RT}$ and $T_{RT}\,$.
\medskip

Now, in a TQFT as in Theorem~\ref{thm-Henn-TQFT} we associate to a full torus
the map ${\bf C}\to{\tt Z}$ whose image is the unit $\,1\,$, and to its complement
in $S^3$ the counit $\varepsilon:{\tt Z}\to{\bf C}$. Hence if an element of the 
mapping class group is represented by $\psi\,:\,{\tt Z}\to{\tt Z}\,$, then
the $\tau_{HKR}$-invariant of the lens-space $L(\psi)$ is 
$\varepsilon\bigl(\psi(1)\bigr)\,$.

In  the relevant part of the presentation of Theorem~\ref{thm-sl2-dec} we find that
the units are given by 
$$
1\;=\;\left[\bar{c}1\\0\ear\right]\!\otimes\,\mu\qquad\quad \varepsilon\;
=\;[1,0]\!\otimes\!\lambda\;,
$$
where $\mu$ is the column vector whose components are the $d(j)$ and $\lambda$
is the projection on the trivial object. (The switch from units to integrals
in the semisimple part occurs since we also switched from a covariant to
a contravariant picture). From the factorized form of the matrices and vectors
we find:

\begin{cor}
Suppose \A=$U_q(s\ell_2)\,$ with $q$ a primitive, (prime) root of unity,
and let $L(\psi)$ denote the lens-space associated to an element 
$\psi\in\pi_0\bigl(Diff(\Sigma_1)^+\bigr)\,$.

Then
$$
\tau_{HKR}\bigl(L(\psi)\bigr)\;=\;\rho\bigl(L(\psi)\bigr)\cdot\tau_{RT}\bigl(L(\psi)\bigr)\quad.
$$
Here, $\rho(M)$ is the order of $H_1(M,{\bf Z})\,$ if that is finite, 
and zero else wise.
\end{cor}

The basic form of this  result is consistent with computations in [KR]
for $L(k,1)$ and $q$ an $8$-th root of unity, especially the linear dependence
on $k=\rho(M)\,$. In particular, we see here  that, even for fixed  finite \A,
 the invariant $\tau_{HKR}\,$ is unbounded. This is very different from the
known semisimple invariants that are all bounded. 
The result also suggests that $\tau_{HKR}\,$ might be computable from
homological data and its semisimple counterpart also in more general 
situations.
We shall leave this as open problem.  

\newpage
 

\setcounter{chapter}{4}

\subsection*{4) Open Questions, and Relations with Cellular Invariants}
\lll{pg-4}
\

In the last chapters we tried to give a survey over those quantum-invariants that
are obtained from surgery presentations of 3-manifolds, but use different algebraic 
tools for the construction. A similar pattern of relations can be detected in the
zoo of cellular invariants, i.e., invariants that use presentations of the manifold 
in terms of a simplicial subdivision, or the cell decomposition given by a 
Heegaard diagram. Since the theory is not as developed here as in the surgical
case (and also since the author is not a true expert on the subject) we shall only
pose a number of questions that arise naturally when we compare the two
families. There are also relations of two different types between cellular
and surgical invariants that we shall briefly talk about.
\medskip

The constructions are all in some way connected to an abelian, 
balanced tensor category. Notions of braiding or modularity are not
needed.  In [TV] 6-j-symbols, which are given as the structure 
constants of $\overline{\Am}^{tr}$, where $\A=U_q(s\ell_2)$, are used
to define state sums over colorings of the edges of a simplicial complex,
representing $M$, which turn out to be invariants. 

The basic reason for invariance is the correspondence between
 subdivision-moves and the  Biedenharn-Elliot (or pentagonal) relations. This
observation  has  already been made  earlier in a non-rigorous 
setting by Ponzano and Regge [PR].  
Let us denote in general an invariant that is obtained in this way from 
semisimple, balanced tensor category, \C,   by $\tau_{TV}\bigl(M,\C\bigr)\,$. 

In [DW] Dijkgraaf and Witten formulated a lattice version of a 
topological field theory associated to a  finite gauge group,  $G$, and  a 3-cocycle,
$\alpha$,  of $G$, such that  the partition function is the order of $Hom(\pi_1(M),G)$,
when $\alpha=1\,$. This is in fact another version of $\tau_{TV}\,$,
if we set  $\C\,=\,C(G)-mod$, with a possibly non-trivial associativity constraint given by
$\alpha\,$, (which is in fact a special 6-j-symbol).

Generalizations of this construction to a more  categorical 
 - but still semisimple - context have been undertaken by many people since.
See, e.g.,  [D...].
\medskip

Another approach is due to G. Kuperberg, [Ku], who associates to the handle attachments
of a Heegaard diagram the multiplication and comultiplication morphisms of a 
Hopf algebra, \A, which again does not have to be  quasi-triangular. The invariant
$\tau_{Ku}(M,\A)\,$ can also be defined on a  ``combed'' manifold, $M$,
if \A is not semisimple (and even without
balancing if we consider framed manifolds, where the notion of balancing in [Ku] differs from the
one in the braided case).

A relation of the form 
$$
\tau_{TV}\bigl(M,\Am\bigr)\;=\;\tau_{Ku}\bigl(M,\A)\qquad{\rm if }\quad\A\quad{\rm semisimple}
$$
has been proven by [BW]. This and also the general ingredients in the constructions 
suggest that $\tau_{TV}$ and $\tau_{Ku}\,$ relate to each other in quite the 
same way as  $\tau_{RT}\,$ and $\tau_{HKR}\,$ do. 
Our first question is thus, whether there is a counterpart to $\tau_L$:

\bQ
Is it possible to define for an abelian, balanced tensor category, \C, 
(that is not necessarily semisimple or braided) an invariant $\tau_X(M,\C)\,$,
which is constructed from a cell-decomposition of $M$, and for which we 
have the specializations $\tau_{TV}=\tau_X$
if \C\ is semisimple, and $\tau_{Ku}=\tau_X$  if \C=\Am\ {\bf ?}
\eQ

Topological quantum field theories for $\tau_{TV}\,$ exists by quite general
arguments. They are constructed by first associating vector spaces of 
``admissible'' colorings  of a triangulation of a  surfaces $\Sigma\,$. We then use 
the state sums to define a ``pre-TQFT'', which respects compositions but not
necessarily the identities.
Reduction of the vector spaces by the projections associated to the cylinders $\Sigma\times[0,1]$
then yields an  actual TQFT. Extending the construction to more general simplicial
subdivisions, it is shown in [KS] that a basis of the vector space 
associated to a surface is
 given by admissible colorings
of the pair of  spines that belong  to the two handle bodies bounding 
$\Sigma\subset S^3\,$. Using the coend-object $F$ in the semisimple setting as in
Lemma~\ref{lm-coend-ss} the result may be expressed in the following form:
\beq\lll{eq-vs-TV}
{\cal V}_{TV}\bigl(\Sigma_g\bigr)\;\cong\;Inv_{\cal C}(F^{\otimes g}\bigr)\otimes 
Inv_{\cal C}(F^{\otimes g}\bigr)^*\;\cong\;End_{\bf C}\bigl(Inv_{\cal C}(F^{\otimes g})\bigr)\qquad\quad.
\eeq

This formula makes sense in any tensor category with limits, even if it is not
semisimple and has no braiding. Although their possibility has been mentioned, 
a construction of TQFT's in the formalism of [Ku] is still missing. Moreover,
we have $\tau_{Ku}\bigl(S^2\times S^1\bigr)\,=\,0\,$ if \A\ is not semisimple,
for the same reasons why $\tau_{HKR}$ vanishes.
The next question is basically about filling this gap for the cellular situation
as it had been done for the surgical in [KL]:

\bQ
Is there a (extended) TQFT,
${\cal V}\,:\,Cob^{conn}(*)_3\,\to\,\vc \;\bigl({\rm or\ }\to \ {\bf AbCat}\bigr)\,$, that 
specializes to $\tau_X$ for closed manifolds   
(or at least to $\tau_X\,=\,\tau_{Ku}\,$) {\bf ?}
\eQ

It has been brought forward in [T] and [Wa] 
that the $\tau_{TV}$ invariant is the norm-square
of $\tau_{RT}\,$ in the case of semisimple categories with hermitian
structures. In fact, an argument of J. Roberts in [R] is easily generalized 
to semisimple categories (without *-structures) to give the identity:
\beq\lll{eq-TV=RT2}
\tau_{TV}\bigl(M,\,\C\bigr)\;=\;\tau_{RT}\bigl(M\#(-M),\,\C\bigr)\;=\;
\tau_{RT}\bigl(M,\,\C\bigr)\tau_{RT}\bigl(-M,\,\C\bigr)\qquad\,.
\eeq

The basic strategy is to replace the 3-dim handle attachments, given by a 
Heegaard diagram, by 4-dim handle attachments of same index, and identify the 
connected sum with the boundary of the resulting four-fold. The application of 
$\tau_{TV}\,$ to the Heegaard diagram is via a ``chain-mail invariant''
shown to be the same as $\tau_{RT}\,$ for the substituted link.

Notice that we need a braided category in order to determine $\tau_{RT}$,
which is superfluous data for the computation of the left hand side of
the equation. In any case we need that $\C$ is semisimple. Yet, it seems
that this property is needed in the proofs only for technical reasons.

\bQ Is is possible to prove  a relation between 
$\tau_X$ and $\tau_{L}\,$ (or at least between $\tau_{Ku}$ and $\tau_{HKR}$) 
that generalizes the one from (\ref{eq-TV=RT2})\ {\bf ?}
\eQ

Notice that the $Inv_{\cal C}(F^g)\,$-spaces that appear in (\ref{eq-vs-TV})
are those associated by the surgery TQFT to a surface. In particular,
 we can rewrite it as 
${\cal V}_{TV}\bigl(\Sigma\bigr)\;=
\;{\cal V}_{RT}\bigl(\Sigma\amalg(-\Sigma)\bigr)$. There is also a natural
identification of the two bases since the ${\cal V}_{RT}$-vector spaces can  also be 
constructed from colorings of  spines.

This situation for vector spaces
has obvious similarities with (\ref{eq-TV=RT2}), (we may 
replace $\#$ by $\amalg$), and the squaring of dimensions is consistent with 
the relation of the invariants for $S^1\times \Sigma\,$. For these reasons
it seems to be likely that there is a TQFT-generalization of 
(\ref{eq-TV=RT2}).

\bQ Is there a  topological ``doubling''-functor ${\bf D}\,:\;Cob_3\,\to\,Cob_3\,$,
which maps a surface $\Sigma$ to $\Sigma\amalg -\Sigma\,$ and a closed manifold 
$M$ to $M\amalg -M\,$ such that:
$$
{\cal V}_{RT}\circ{\bf D}\;=\;{\cal V}_{TV}
$$
for a fixed braided category \C\ {\bf ?} What about  non-semisimple 
generalizations {\bf ?}
\eQ

Constructive answers to the above questions would give us  a way to compute
all known quantum-invariants and the corresponding TQFT's 
 from $\tau_L$, if we consider only cellular invariants that start from braided
categories. Finding proofs should only be a matter of  gaining a thorough 
and sufficiently deep understanding of the existing constructions.
 
There is, however, also an anticipated  relation between cellular and surgical
invariants that would allow us also to compute the cellular ones for 
non-braided categories from $\tau_L\,$. Unlike the previous relation
there are still a few conceptual ideas missing to make it work in full
generality.

The relation is with respect to a fixed manifold
$M$. Instead we change the algebraic data in a non-trivial way.
Specifically, if $\tau_{Ku}\,$ is defined for a Hopf algebra $\cal B$
we can define $\tau_{HKR}\,$ for the corresponding Drinfel'd 
double $D({\cal B})$, since this is quasi-triangular and, moreover, 
always modular. The notion of  doubles has been extended by Majid [M2]
to  balanced,
monoidal categories \C, such that $D(\C)\,$ is a balanced, modular,
braided tensor category. Following a conjecture of Gelfand and Kazhdan,
we thus ask the next question.

\bQ Is there a relation of the form
$$
\tau_L\bigl(M,D({\cal C})\bigr)\;=\;\tau_X\bigl(M,\C\bigr)
$$
for general balanced, abelian tensor categories {\bf ?} 
What about TQFT's {\bf ?}
\eQ

There is in fact more than just aesthetical evidence for this fact:

First of all we already know that this relation is true for finite groups,
comparing the results from [DW] for $\tau_{TV}(G)\,$ and from
[AC] for $\tau_{RT}(D(G))\,$, where the answer is in both cases the order of
$Hom\bigl(\pi_1(M),G\bigr)\,$. The vector spaces of 
${\cal V}_{TV}\bigl(G\bigr)\,$ and 
${\cal V}_{RT}\bigl(D(G)\bigr)\,$ are identified by the easy but remarkable 
correspondence  
$
C\Bigl(Hom\bigl(\pi_1(\Sigma_g),G\bigr)/G\Bigr)\,\cong\,\Bigl(D(G)^{\otimes}\Bigr)_{ad}\,  
$, which is given on representatives by 
 
$$
\phi \qquad\mapsto\quad\sum_{\gamma\in Hom\bigl(\pi_1(\Sigma_g),G\bigr)}
\phi(\gamma)\;
\gamma(a_1)\delta_{\gamma(b_1)}\otimes\ldots\otimes \gamma(a_g)\delta_{\gamma(b_g)}
$$

Here, the $a_j$ and $b_j$ are the usual generators of $\pi_1(\Sigma)\,$, and 
the basis vector $\delta_a\in C(G)\,$ is the characteristic function of 
$\{a\}$.
\medskip

In the case of a general algebra, \A, we expect for ${\cal V}_{Ku}\,$ 
to associate an algebra, \A,
to a generating $a$-cycle of $\Sigma$ and the dual $\A^*$ if the cycle lives
in an opposite Lagrangian subspace. (Let us omit the passage to invariances,
i.e., imagine we have introduced a puncture.)  This would then reproduce the
${\cal V}_{HKR}\,$ vector spaces. 

This also fits in nicely with the actions of the action of the $S^*$ and $T^*$
matrices for doubles, from (\ref{eq-ST-1}). As the $S$-transformation flips 
$a$-cycles into $b$-cycles we expect it to be represented by a permutation and
other canonical isomorphisms. From the ${\cal V}_L$-picture this is 
confirmed by the result in (\ref{eq-S-fact}).
\medskip

The relation can obviously be combined with the previous one,
giving rise to a formula for $\tau_L\,$ only, with \C\  a braided category:
\beq\lll{eq-2M=DC}
\tau_L\bigl(M,\,D(\C)\bigr)\quad=\quad\tau_L\bigl(M\amalg-M,\,\C\bigr)
\eeq

This would imply interesting formulae for quasi-triangular Hopf algebras.
For example the dimension of vector spaces would  obey 
$dim\Bigl({\tt Z}_{D({\cal A})}\Bigr)\,=\,dim\Bigl({\tt Z}_{\cal A}\Bigr)^2\,$,
which is wrong for plenty of examples of Hopf algebras, \A, that are not 
quasi-triangular.
\medskip

More abstractly, if we think of a TQFT  as a functor, 
$Cob_3\times{\bf AbCat}\,\to\,{\bf AbCat}\;$, the relation in (\ref{eq-2M=DC})
suggests that  there are topological operations, like the ``doubling''-functor
from Question {\em 4}, which are dual to the doubling of abelian categories.
Let us thus conclude  with the following  esoterical  question:

\bQ Is it possible to describe a topological doubling in a way that 
it is dual to the doubling of abelian categories, and consistent with
an extended ${\cal V}_L\,$ {\bf ?}
\eQ

\medskip

\subsection*{Acknowledgements:} I thank  the organizers of
the Aarhus meeting for the opportunity to give a talk. Also, I wish to thank
all those mentioned in this survey, who had the patience to explain 
their constructions of three-manifold invariants 
 to me; especially, V. Lyubashenko, but also L. Kauffman, D. Kazhdan, 
G. Kuperberg, D. Radford, and V. Turaev. The author is supported by 
NSF-grant DMS 9304580.

\medskip


\subsection*{References}
{\small
\ben

\item[[A]] Atiyah, M.: Publ. Math. Inst. Hautes Etudes Sci. {\bf 68}, 175-186 (1989).

\item[[AC]] Altschuler, D., Coste, A.: Commun. Math. Phys. {\bf 150}, 83-107 (1992).
{\tt hep-th/9202047}.

\item[[BW]] Barrett, J., Westbury, B.: Math. Proc. Camb. Phil. Soc. {\bf 118} (1995).

\item[[CY]] Crane, L., Yetter, D.: On Algebraic Structures Implicit in
Topological Quantum Field Theories, Preprint {\tt hep-th 9412025}.

\item[[D]] Deligne, P.: In ``The Grothendieck Festschrift''  {\bf II}, Progr. in Math., Birkh\"auser (1990).

\item[[Dr]] Drinfel'd, V.: Leningrad Math. J., {\bf 1} No.2, 321-342 (1990).

\item[[D...]] Durhuus, B., Jacobson, H., Nest, R.: Reviews in Math. Physics {\bf 5}, 1-67 (1993). 

Felder, G., Grandjean, O.: On Combinatorial 3-Manifold Invariants, Preprint (1992). 

Gelfand, S., Kazhdan, D.: Invariants of 3-Dimensional Manifolds, Preprint (1994).

Yetter, D.: Top. and its App. {\bf 58}, 47-80 (1994).

Barrett, J., Westbury, B.: Invariants of Piecewise Linear Manifolds. {\tt hep-th/9311155}.

\item[[DPR]] Dijkgraaf, R., Pasquier, V., Roche, P.: Nucl. Phys. B. Proc. Suppl. {\bf 18B}, 60-72 (1990).

\item[[DW]] Dijkgraaf, R.,  Witten, E.: Commun. Math. Phys. {\bf 129}, 393-429 (1990).

\item[[F]] Freed, D.: Commun. Math. Phys. {\bf 159}, 343-398 (1994).
\newline $\;$ Freed, D, Quinn, F.: Chern Simons Theory with finite Gauge Group, Preprint (1992).

\item[[FLM]] Lyubashenko, V., Majid, S.: J. Algebra {\bf 166}, 506-528 (1994).

Crivelli, M., Felder, G., Wieczerkowski, C.:  Lett. Math. Phys. {\bf 30} (1994).

\item[[H]] Hennings, M.: Invariants of Links and 3-Manifolds obtained from Hopf Algebras,
 Preprint (1990).

\item[[Ka]] Kazhdan, D.: Private communication, Harvard lecture notes.

\item[[Ke1]] Kerler, T.: In "New Symmetry Principles in QFT", Physics  {\bf 295}, Plenum Press (1991).

\item[[Ke2]] Kerler, T.: Commun. Math. Phys. {\bf 168},  353-388 (1994).

\item[[Ke3]] Kerler, T.: On the Connectivity of Cobordisms, and
Half-Projective TQFT's , Preprint (1995),
submitted to Contemp. Math.

\item[[Ke4]]  Kerler, T.: Hopf algebras in Cobordisms-Categories, 
in preparation. 

\item[[Ke5]]  Kerler, T.: Bridged Links and Tangle Presentations of
Cobordism Categories, Preprint (1994).

\item[[KL]] Kerler,T., Lyubashenko, V.: Non-semisimple TQFT's for Connected Surfaces, Preprint (1995).

\item[[KR]] Kauffman, L.,  Radford, D.: Invariants of 3-Manifolds Derived from Finite
Dimensional Hopf Algebras, Preprint {\tt hep-th/9406065}; 
to appear in J. Knot Theory and its Rami.

\item[[Ki]] Kirby, R.: Invent. Math. {\bf 45}, 35-56 (1978).

\item[[KS]] Karowski, M., Schrader, R.: Commun. Math. Phys. {\bf 151}, 355-402 (1993).

\item[[Ku]] Kuperberg, G.: Non-Involutory Hopf algebras \& 3-Manifold Invariants, Preprint (1994).

\item[[L]] Lyubashenko, V.: Commun. Math. Phys.  {\bf 172}, 467-516 (1995).

Lyubashenko, V.: J. Pure Appl. Alg. {\bf 98}, 279-327 (1995)

\item[[LS]] Larson, R., Sweedler, M.: Amer. J. of Math., {\bf 91}, No 1., 75-94 (1969).

\item[[McL]] MacLane, S.:``Categories for the Working Mathematician", Springer-Verlag (1971).

\item[[M1]]  Majid, S.: J. Pure Appl. Algebra {\bf 86}, 187-221 (1993).

\item[[M2]]  Majid, S.: Rend. Circ. Mat. Palermo Suppl. {\bf 26},197-206 (1991).

\item[[O]] Ohtsuki, T.: Math. Proc. Cambridge Philos. Soc. {\bf 117}, 
259-273 (1995).

\item[[PR]] Ponzano, G., Regge, T.: In ``Spectroscopic and Group Theoretical Methods in
Physics'',  North-Holland, 1-58 (1968). 

\item [[Rd]] Radford, D.: Amer. J. Math. {\bf 98}, 333-355 (1976).

\item[[Re]]Reshetikhin, N.: Leningrad Math. J. {\bf 1}, No. 2,  491-513 (1990).

\item[[RT]] Reshetikhin, N., Turaev, V.: Invent. Math. {\bf 103}, 547-597 (1991).

\item [[R]]Roberts, J.: Skein Theory and Turaev-Viro Invariants, Preprint (1993). 

\item[[Sw]] Sweedler, M.: Annals of Math. {\bf 89}, 323-335 (1969).

\item[[T]] Turaev, V.: ``Quantum Invariants of 3-Manifolds''. Walter de Gryter (1994). 

\item[[Ta]] Saavedra Rivano, N.: Categories Tannakiennes, LNM {\bf 265}, 
Springer (1972).

Ulbricht, K.-H.: Israel J. Math. {\bf 72}, 252-256 (1990).

Majid, S.: Internat. J. Mod. Phys. {\bf 6}, 4359-4374 (1991).

\item[[TV]] Turaev, V., Viro, O.: Topology {\bf 31}, 865-902 (1992).

\item[[Wa]] Walker, K.: On Witten's 3-Manifold Invariants, Preprint (1991).
 
\item[[Wi]] Witten, E.: Commun. Math. Phys. {\bf 121}, 351-399 (1989).

\een
\bigskip
}
{\small {\sc Institute for Advanced Study, Olden Lane,
Princeton,  NJ 08540, USA}

{\em E-mail-address:} kerler@math.ias.edu }

\end{document}